\input amstex
\documentstyle{amsppt}
\voffset=-3pc
\def\bC{\Bbb C}
\def\bM{\Bbb M}
\def\bR{\Bbb R}

\def\bK{\Bbb K}

\def\cD{\Cal D}
\def\cB{\Cal B}
\def\cP{\Cal P}

\def\ob{\overline b}

\def\ot{\overline t}

\def\al{\alpha}
\def\del{\delta}
\def\ep{\epsilon}
\def\lam{\lambda}
\def\var{\varphi}
\def\ga{\gamma}
\def\sa{{\text{sa}}}
\def\msa{M_{\text{sa}}}

\def\ssa{S_{\text{sa}}}
\def\asa{A_{\text{sa}}}

\def\m{{\text{m}}}

\def\tA{\widetilde A}
\def\tx{\widetilde x}
\def\ty{\widetilde y}
\def\tb{\widetilde b}
\def\ta{\widetilde a}
\def\tc{\widetilde c}
\def\tf{\widetilde f}
\def\th{\widetilde h}
\def\tk{\widetilde k}
\def\tp{\widetilde p}
\def\Asa{\widetilde A_{\text{sa}}}

\def\fp{F(p)}
\def\b1{\bold 1}

\magnification=\magstep1
\parskip=6pt
\NoBlackBoxes
\topmatter
\title Semicontinuity and Closed Faces of C*--Algebras
\endtitle  
\author Lawrence G.~Brown
\endauthor
\keywords{Operator algebras, Semicontinuity, Closed projection, Operator convex}\endkeywords
\subjclass{Primary 46L05}\endsubjclass

\abstract{C. Akemann and G. Pedersen defined three concepts of semicontinuity for self--adjoint elements of $A^{**}$, the enveloping von Neumann algebra of a $C^*$--algebra $A$. We give the basic properties of the analogous concepts for elements of $pA^{**}p$, where $p$ is a closed projection in $A^{**}$. In other words, in place of affine functionals on $Q$, the quasi--state space of $A$, we consider functionals on $F(p)$, the closed face of $Q$ suppported by $p$. We prove an interpolation theorem: If $h\geq k$, where $h$ is lower semicontinuous on $F(p)$ and $k$ upper semicontinuous, then there is a continuous affine functional $x$ on $F(p)$ such that $k\leq x\leq h$. We also prove an interpolation--extension theorem: Now $h$ and $k$ are given on $Q$, $x$ is given on $F(p)$ between $h_{|\fp}$ and $k_{|\fp}$, and we seek to extend $x$ to $\tx$ on $Q$ so that $k\leq\tx\leq h$. We give a characterization of $pM(A)_{\sa}p$ in terms of semicontinuity. And we give new characterizations of operator convexity and strong operator convexity in terms of semicontinuity.}
\endabstract
\endtopmatter

\S1. Definitions, notations, and basic properties.

\noindent
For a $C^*$--algebra $A$, $S=S(A)$ denotes the state space of $A$ and $Q=Q(A)$ the quasi--state space, $Q(A)=\{\var\in A^*:\var\geq 0 \;\text{and}\;\|\var\|\leq 1\}$. 
E. Effros \cite{E} showed that norm closed faces of $Q(A)$ containing 0 correspond one--to--one to projections $p$ in $A^{**}$, via $\fp=\{\var\in Q: \var(\b1-p)=0\}$. 
Then $p$ is called {\it closed} if $\fp$ is weak$^*$ closed and {\it open} if $p$ is the support projection of a hereditary $C^*$--subalgebra of $A$. 
It was proved by Effros in \cite{E} that these definitions imply $p$ is closed if and only if $1-p$ is open, and the open/closed terminology was introduced by Akemann in \cite{A1}. 
Akemann \cite{A2} introduced the concept of compact projections. 
The projection $p$ is {\it compact} if and only if $\fp\cap S$ is weak$^*$ closed. 
Among several equivalences, $p$ is compact if and only if it is closed and there is $a$ in $A$ such that $p \leq a\leq\b1$; and $p$ is compact if and only if it is closed in ${\tA}^{**}$. 
Here $\tA=A+\bC\b1$, where $\b1$ is the identity of $A^{**}$. 
For the semicontinuity concepts introduced by Akemann and Pedersen in \cite{AP} we need some more notations. 
For $S\subset A^{**}$, $\ssa=\{h\in S: h^*=h\}, S_+=\{h\in S: h\geq 0\}$, $S^{\m}$ is the set of limits in $A^{**}$ of bounded increasing nets from $S$, $S_{\m}$ the set of limits of bounded decreasing nets from $S$, and $S^\sigma$ and $S_\sigma$ are defined similarly using monotone sequences. 
Also $^-$ denotes norm closure. 
Then $h$ is {\it strongly}  {\it lower semicontinuous} ({\it{lsc}}) if $h\in((\asa)^{\m})^-$, {\it middle}  {\it lsc} if $h\in(\Asa)^{\m}$, and {\it weakly} {\it lsc} if $h\in((\Asa)^{\m})^-$. 
The upper semicontinuous (usc) concepts are defined analogously, using ${}_{\m}$ instead of ${}^{\m}$, and $h$ is usc in any sense if and only if $-h$ is lsc in the same sense. 
We use completely analogous definitions for semicontinuity on $\fp$. 
If $p$ is a closed projection in $A^{**}$ and $h$ is in $pA^{**}_{\sa}p$, then $h$ is {\it strongly} {\it lsc} {\it on} $p$ if $h\in((p\asa p)^{\m})^-$, $h$ is {\it middle} {\it lsc} {\it on} $p$ if $h\in(p\Asa p)^{\m}$, and $h$ is {\it weakly} {\it lsc} {\it on} $p$ if $h\in((p\Asa p)^{\m})^-$; and the usc concepts are defined similarly.

Under the Kadison function representation, $A^{**}$ is identified with the space of bounded affine functionals vanishing at 0 on $Q$ and $A$ is identified with the subspace of weak$^*$ continuous functionals. These identifications are isometric on $\asa^{**}$ and $\asa$. Since $\fp$ is the normal quasi--state space of the von Neumann algebra $pA^{**}p$, similar identifications are available relative to a closed projection $p$. Thus $p\asa^{**}p$ is identified with $B_0(p)$, the space of $\bR$--valued bounded affine functionals on $\fp$, vanishing at 0, and $p\asa^{**}p$ is identified with $A_0(p)=\{f\in B_0(p):f\; \text{is weak$^{*}$ continuous}\}$. (An argument justifying the last fact, which is presumably folklore, was provided in the proof of [B1, Corollary 3.5].) Usually these identifications will be used without any notation, but when a notation is necessary, $\hat a$ will denote the functional on $Q$ corresponding to $a$ in $A^{**}$, and the unwieldy symbol $F_p(a)$ will denote the functional on $\fp$ corresponding to $a$ in $pA^{**}p$. Note that $F_p(pap)=\hat a_{|\fp}$.

Since $F_p(p)(\var)=\|\var\|$, for $\var\in\fp, F_p(p)$ is an lsc function on $\fp$. If $p$ is compact, then by the above, $p\in pAp$ and hence $F_p(p)$ is continuous. Otherwise, $\fp\cap S$ is not weak$^*$ closed. Thus if $p$ is not compact, $F_p(p)$ is not usc and $p$ is not in $pAp$. Once it is proved that strong lsc on $p$ implies middle lsc on $p$, it will then be clear that if $p$ is compact, all three types of semicontinuity coincide, and if $p$ is not compact, the strong and middle types definitely differ. 

Some of the Akemann--Pedersen work is part of the abstract theory of compact convex sets. For example, see the beginning of [AP,\S 3]. Thus in many cases, their arguments apply in our situation and our proofs can be abridged. In this connection note that the functional $F_p(p)$, which plays the same role here that $\hat\b1$ plays in \cite{AP}, is determined by the structure of $\fp$ as a convex set and the choice of the complemented extreme point 0.

\proclaim{Lemma 1.1} (cf. [AP, Theorem 2.1]) If $p$ is a closed projection in $A^{**}$ and $h$ is in $p\asa^{**}p$ then $F_p(h)$ is lsc, if and only if $h$ is the $\sigma$--weak limit of an increasing net $(h_i)$, where $h_i=\lam_ip+pa_ip\in p\Asa p$ and $\lim \lam_i=0$. In this case $(\lam_i)$ may be taken strictly increasing.
\endproclaim 

\demo{Proof}The argument from \cite{AP} applies unchanged.
\enddemo

We state a known result for reference.

\proclaim{Lemma 1.2} ([B1, Corollary 3.4]) If $p$ is a closed projection in $A^{**}$, $h$ is in $p\asa p$, and $-sp\leq h\leq tp$ for $s,t\geq 0$, then there is $a$ in $\asa$ such that $pap=h$ and $-s\b1\leq a\leq t\b1$. (If $A$ is unital, $s$ and $t$ may be arbitrary.)
\endproclaim 

\proclaim{Proposition 1.3} (cf. [AP, Theorem 2.1 and Proposition 2.2]) If $p$ is a closed projection in $A^{**}$ and $h$ is on $p\asa^{**}p$, then the following are equivalent:

(i) $h$ induces an lsc function on $\fp$.

(ii) $h$ the is $\sigma$--weak limit of an increasing net $(h_i)$ such that $h_i=\lam_i p+pa_ip\in p\Asa p$, and $\lim \lam_i=0$.

(iii) $\forall \ep>0,\, h+\ep p \in (p\asa p)^{\m}$. 

(iv) $h\in((p\asa p)^{\m})^-$.

\noindent
Moreover, if $h\geq 0$ and $h$ is strongly lsc on $p$, then $\forall\ep>0,\, h+\ep p \in (pA_{+p})^{\m}$.
\endproclaim

\demo{Proof}The proof is almost the same as that in \cite{AP}, but one additional argument is needed for the proof that (ii) implies (iii) and the proof of the last sentence. Let $\cP=\{a\in A_+:\|a\|<1\}$, the canonical approximate identity of $A$ ([P2, p.11]). If $(h_i)$ is as in (ii), we construct a net $k_{i,b}=(\lam_i+\ep)b+pa_ip$, $b\in p\cP p$, where we use only the pairs  $(i,b)$ with $\lam_i+\ep>0$. For the last sentence we also require that $k_{i,b}\geq 0$. The ordering on the index set is modified as follows: $(i,b)\leq(j,c)$ if and only if $i\leq j$ and $k_{i,b}\leq k_{j,c}$. In order to make the Akemann--Pedersen argument work, we need the following: If $x\in p\asa p$ and $\sigma_{pA^{**}p}(x)\subset (-\infty,t]$ with $0<t<1$, then there is $b$ in $p\cP p$ such that $x\leq b$. To prove this, first cite 1.2 to write $x=p\tx p$ with $\tx\in\asa$ and $\sigma(\tx)\subset(-\infty, t]$ and then let $b=p\tx_+ p$. One can also cite 1.2 for the notationally convenient fact that $pA_+p=(pAp)_+$. 

When $p$ is closed but not compact and $a\in \asa$, the condition $\lam p+pap\geq 0$ does not necessarily imply $\lam\geq 0$. In other words, it is possible that
$$\exists a\in A_+ \quad \text{such that}\quad p\leq pap.\leqno(1)$$
The condition (1) was studied in the author's unpublished manuscript \cite{B2}, where it was shown (for $p$ closed) to be equivalent to the existence of a compact projection $q$ such that $\|p-q\|<1$. Also $pAp$ is (non-isometrically) completely order isomorphic to $qAq$. The following easy lemma is needed for this situation.
\enddemo

\proclaim{Lemma 1.4} If $p$ is a closed but not compact projection in $A^{**}$, then there is a constant $K$ such that $|\lam|\leq K\|\lam p+pap\|$ for $a\in A$.
\endproclaim 
\demo{Proof}This follows from elementary Banach space theory, the fact that $pAp$ is closed ([APT, Proposition 4.4]), and the fact that $p\notin pAp$.
\enddemo

\proclaim{Proposition 1.5} (cf. [AP, Proposition 2.5]) If $p$ is a closed projection in $A^{**}$ and $h\in p\asa^{**}p$, then $h\in (p\Asa p)^{\m}$ if and only if there is $\lam$ such that $h+\lam p$ is strongly lsc on $p$. 
\endproclaim 
\demo{Proof}If $h+\lam p$ is strongly lsc, then Proposition 1.3 clearly implies that $h\in(p\Asa p)^\m$. Conversely suppose $h$ is the limit of an increasing net $(h_i)$, where $h_i=\lam_i p+pa_ip$. If $p$ is non--compact, then by Lemma 1.4, there is $\lam$ such that $\lam+\lam_i\geq 0$, for $i\geq i_0$. Thus each such $h_i+\lam p$ gives an lsc functional on $\fp$, and so does $h+\lam p$. Then the equivalence of (i) and (iv) in Proposition 1.3 shows that $h+\lam p$ is strongly lsc. The compact case is trivial.
\enddemo

\example{Remark 1.6} It was shown in [AP, Theorem 3.3] that $h$ in $\asa^{**}$ is weakly lsc if and only if it induces an lsc function on $S$. The analogous result for closed faces is false as shown below in Example 3.4. A portion of the Akemann--Pedersen argument is very general and still applies here. Namely, $h$ in $p\asa^{**}p$ is lsc on $\fp\cap S$ if and only if there is a net $(h_i)$ in $p\Asa p$ such that $h_i\leq h$, $\forall i$, and $(h_i)$ converges to $h$ $\sigma$--weakly. (The net may not be bounded.) The rest of the Akemann--Pedersen argument fails since $p\tA p$ may not be an algebra and hence may not be inverse closed.

It may be imagined that we should have defined weak semicontinuity (or that we should define a fourth and even weaker kind of semicontinuity) to mean that $h$ is a semicontinuous function on $\fp\cap S$. But we believe that this last condition is too weak to be useful, except when it implies $h\in((p\Asa p)^{\m})^-$. The following open question seems mildly interesting: Is $\{h\in p\asa^{**}p: h$ is lsc on $\fp\cap S\}$ the smallest subset of $p\asa^{**}p$ containing $((p\Asa p)^{\m})^-$ and closed under increasing convergence? 

Another basic subject is semicontinuity of projections. Akemann and Pedersen showed ([AP, Theorem 3.6]) that $p$ is open if and only if it is strongly lsc if and only if it is weakly lsc, whence $p$ is closed if and only if it is middle usc if and only if it is weakly usc. And the author observed ([B1, 2.47]) that $p$ is compact if and only if it is strongly usc. Only a portion of this goes through in the relative case.
\endexample

\proclaim{Proposition 1.7} If $p$ is a closed projection in $A^{**}$ and $q$ is a subprojection of $p$, then $q$ is strongly lsc on $p$ if and only if $p-q$ is closed (i.e., $q$ is relatively open), and $q$ is strongly usc on $p$ if and only if $q$ is compact. 
\endproclaim 

\demo{Proof}If $q$ is strongly lsc on $p$, then $F(p-q)=\{\var \in \fp:F_p(q)(\var)\leq 0\}$, which is closed. Conversely, if $p-q$ is closed, then $\b1-p+q$ is open and hence strongly lsc in $A^{**}$. {\it{A}} {\it{fortiori}}, $q=p(1-p+q)p$ is strongly lsc on $p$. If $q$ is strongly usc on $p$, then $S\cap F(q)=\{\var\in\fp:F_p(q)(\var)\geq 1\}$, which is closed. Therefore $q$ is compact. Conversely, if $q$ is compact, then $q$ is strongly usc in $A^{**}$ and, {\it {a fortiori}}, strongly usc on $F(p)$.
\enddemo

\example{Example 1.8}Let $A=c\otimes\bK$, the algebra of norm convergent sequences in $\bK$, the set of compact operators on $l^2$. 
Then $A^{**}$ can be identified with the algebra of bounded indexed collections, $(h_n)_{1\leq n\leq\infty}$, with each $h_n$ in $B(l^2)$. 
Let $v_n={1\over\sqrt{2}}e_1+{1\over\sqrt{2}}e_{n+1}$, where $e_1, e_2,\dots$ are the standard basis vectors in $l^2$. 
Let the closed projection $p$ in $A^{**}$ be given by $p_n=v_n\times v_n$ for $n<\infty$, and $p_\infty=e_1\times e_1$. 
Here $v\times w$ denotes the rank one operator $x\mapsto(x,w)v$. 
An element $h$ of $p\asa^{**}p$ is given by a bounded collection $(t_n)_{1\leq n\leq\infty}$ of real numbers such that $h_n=t_np_n$. 
Then it is easily seen that $h$ is strongly lsc on $p$ if and only if ${1\over 2}t_\infty\leq\lim\inf t_n$ and strongly usc on $p$ if and only if ${1\over 2}t_\infty\geq\lim\sup t_n$. 
It follows that every element of $p\asa^{**}p$ is middle lsc and middle usc on $p$. 
Since there are subprojections of $p$ which are not closed, for example $p-p_\infty$, this shows that middle semicontinuity on $p$, for a subprojection $q$ of $p$, does not imply that $q$ is closed or relatively open.
\endexample

\S2. The interpolation and interpolation--extension theorems.

\proclaim{Lemma 2.1} Assume $p$ is a closed projection in $A^{**}$, $h$ is strongly lsc on $p$, $k$ is strongly usc on $p$, and $h\geq k$. Then there is a function $f$ such that $\lim\limits_{\ep\to 0^+}f(\ep)=0$ and: If $\delta>0, x\in p\asa p$, and $k-\ep p\leq x\leq h+\ep p$, then $\exists y\in p\asa p$ such that $\|y-x\|\leq f(\ep)$ and $k-\delta p\leq y\leq h+\delta p$.
\endproclaim 

\demo{Proof}The proof has a lot in common with that of [B1, Lemma 3.14], but there are enough extra steps needed that it seems advisable to write out a complete version. Choose a net $(\lam_\al p+pa_\al p)$ which increases to $h$ and has all the properties of Lemma 1.1. Choose a similar net $(\mu_\beta p+pb_\beta p)$ which decreases to $k$. Let $\del>0$. Then with the help of Dini's theorem we can see that for sufficiently large $\al$ and $\beta$, all of the following hold:
$$\align
(\mu_\beta-\ep-\del)p+pb_\beta p&\leq x\leq(\lam_\al+\ep+\del)p+pa_\al p,\\
\mu_\beta p+pb_\beta p&\leq(\lam_\al+\del)p+pa_\al p,\quad\text{and}\\
\mu_\beta-\lam_\al&\leq\del.
\endalign
$$
Fix one such $(\al,\beta)$. Since $\lam_\al<0$ and $\mu_\beta>0$, we also have:
$$\align
pb_\beta p-(\ep+\del)p&\leq x\leq pa_\al p+(\ep+\del)p,\\
pb_\beta p&\leq pa_\al p+\del p,\\
p(a_\al-b_\beta)p&\leq h-\lam_\al p-k+\mu_\beta p\leq h-k+\del p,\\
pa_\al p&\leq h+\del p, \quad\text{and}\\
pb_\beta p&\geq k-\del p.
\endalign
$$

By 1.2 there is $\tc$ in $\asa$ such that $p\tc p=p(a_\al-b_\beta)p$ and $-\del\b1\leq\tc\leq(\|h-k\|+\del)\b1$. Also let $\tb=b_\beta$, $\ta=\tb+\tc$, and choose $\tx$ in $\asa$ so that $p\tx p=x$. Then 
$$p(\tb-\ep\b1-\del\b1)p\leq p\tx p\leq p(\ta+\ep\b1+\del\b1)p,$$
and hence $\exists\lam>0$ such that
$$\align
-\lam(\b1-p)+\tb-(\ep+2\del)\b1&\leq\tx\leq\ta+(\ep+2\del)\b1+\lam(\b1-p),\quad\text{and}\\
-\lam(\b1-p) +\tb-\del\b1&\leq\ta +(1/2)\del\b1+\lam(\b1-p).
\endalign
$$
Let $(e_\ga)$ be an approximate identity of her$(\b1-p)$, the hereditary $C^*$--subalgebra of $A$ supported by $\b1-p$. Then by Dini's theorem, for $\ga$ sufficiently large we have:
$$\align
-\lam e_\ga+\tb-(\ep+3\del)\b1&\leq\tx\leq\ta+(\ep+3\del)\b1+\lam e_\ga,\quad\text{and}\\
-\lam e_\ga+\tb-\del\b1 &\leq\ta +\del\b1 +\lam e_\ga.
\endalign
$$
Fix one such $\ga$. Then
$$\tx=-\lam e_\ga+\tb-(\ep+3\del)\b1+(\ta-\tb+(2\ep+6\del)\b1+2\lam e_\ga)^{1\over 2}t(\ta-\tb+(2\ep+6\del)\b1+2\lam e_\ga)^{1\over 2},$$
where $0\leq t\leq\b1$ and $t\in{1\over 2}\b1+A\subset\tA$.\newline
Let
$$\ty=-\lam e_\ga+\tb-\del\b1+(\ta-\tb+2\del\b1+2\lam e_\ga)^{1\over 2}t(\ta-\tb+2\del\b1+2\lam e_\ga)^{1\over 2},$$
and $y=p\ty p$. Then $\ty\in\asa$ and
$$pb_\beta p-\del p\leq y\leq p a_\al p+\del p.$$
Hence
$$k-2\del p\leq y\leq h+2\del p.$$
Also,
$$\|y-x\|\leq\ep+2\del+\sqrt{2\ep+4\del}\bigg[\|(\ta-\tb+(2\ep+6\del)\b1+2\lam e_\ga)^{1\over 2}p\|+\|p(\ta-\tb+2\del\b1+2\lam e_\ga)^{1\over 2}\|\bigg].$$
Since $\|z^{1\over 2}p\|=\|pzp\|^{1\over 2}$ and $\|pz^{1\over  2}\|=\|pzp\|^{1\over 2}$,
$$\align
\|y-x\|&\leq\ep+2\del+\sqrt{2\ep+4\del}\bigg[(\|h-k\|+2\ep+7\del)^{1\over 2}+(\|h-k\|+3\del)^{1\over 2}\bigg]\\
&\leq\ep+2\del+\sqrt{2\ep+4\del}\bigg[2\|h-k\|^{1\over 2}+(2\ep+7\del)^{1\over 2}+(3\del)^{1\over 2}\bigg]\\
&\leq\max(C_1\ep, C_2\ep^{1\over 2}\|h-k\|^{1\over 2}),
\endalign
$$
if $\del$ is small enough.
\enddemo

\proclaim{Theorem 2.2}Assume $p$ is a closed projection in $A^{**}$, $h$ is strongly lsc on $p$, and $k$ is strongly usc on $p$, and $k\leq h$. 

(i) Then there is $x$ in $p\asa p$ such that $k\leq x\leq h$.

(ii) Moreover, there is a function $f$ such that $\lim\limits_{\ep\to 0^+} f(\ep)=0$ and for each $y$ in $p\asa p$ such that $k-\ep p\leq y\leq h+\ep p$, there is $x$ in $p\asa p$ with $\|x-y\|\leq f(\ep)$ and $k\leq x\leq h$.

(iii) If $\Cal S=\{x\in p\asa p:k\leq x\leq h\}$ and $\Cal T=\{y\in p\asa^{**}p:k\leq y\leq h\}$, then $\Cal S$ is $\sigma$--weakly dense in $\Cal T$.
\endproclaim 

\demo{Proof} The deduction of parts (i) and (ii) from the lemma is routine and is identical to the deduction of [B1, Theorem 3.15] from [B1, Lemma 3.14]. Note that the initial $\ep$ need not be small.

(iii) Let $y$ be in $\Cal T$ and let $V$ be a symmetric, convex, $\sigma$--weak neighborhood of 0. We need to find $x$ in $\Cal S$ such that $x-y\in V$. First choose $\del>0$ such that the ball of radius $f(3\del)$ is contained in $1/3\, V$, for the $f$ of part (ii). 
Choose nets $(\lam_\al p+pa_\al p)$ and $(\mu_\beta p+pb_\beta p)$ as in the proof of the lemma. As before, for sufficiently large $\al$ and $\beta$, we have
$$\align
pb_\beta p&\leq pa_\al p+\del p,\\
(2)\qquad\qquad pa_\al p&\leq h+\del p,\quad\text{and}\\
pb_\beta p&\geq k-\del p.
\endalign
$$
Also note that $(pa_\al p)$ and $(pb_\beta p)$ converge $\sigma$--strongly to $h$ and $k$ and that $k-2\del p\leq y\leq h+2\del p$.

There is $t$ in $p\asa^{**}p$ such that 
$$0\leq t\leq p\quad\text{and}\quad y=k-2\del p+(h-k+4\del p)^{1\over 2}t(h-k+4\del p)^{1\over 2}.$$
Let $z_{\al\beta}=(pb_\beta p)-2\del p+(pa_\al p-pb_\beta p+4\del p)^{1\over 2}t(pa_\al p-pb_\beta p+4\del p)^{1\over 2}$, for $\al,\beta$ sufficiently large. Since $z_{\al\beta}\to y$, we may fix $(\al,\beta)$ so that (2) holds and $z_{\al\beta}- y\in 1/3\,V$. Then
$$
pb_\beta p-2\del p\leq z_{\al\beta}\leq pa_\al p+2\del p.$$

Now, as above, choose $\tc$ in $\asa$ such that $p\tc p=p(a_\al-b_\beta)p$ and $\tc\geq -\del \b1$, let $\tb=b_\beta$, and let $\ta=\tb+\tc$. Then choose $\lam>0$ so that
$$
\tb-\del\b1\leq\del\b1+\ta+\lam(\b1-p).$$
Finally, choose $e_\ga$ from an approximite identity of her$(\b1-p)$ so that
$$\tb-2\del\b1\leq 2\del\b1+\ta+\lam e_\ga.$$

Then there is $s$ in $\asa^{**}$ such that $0\leq s\leq \b1$ and 
$$
z_{\al\beta}=pb_\beta p-2\del p+p(\ta-\tb+4\del\b1+\lam e_\ga)^{1\over 2}s(\ta-\tb+4\del\b1+\lam e_\ga)^{1\over 2}p.$$
By the Kaplasky density theorem, there is a net $(r_i)$ in $\asa$ such that $\|r_i\|\leq 1/2$ and $r_i\to s-{1/2}\,\b1$, $\sigma$--strongly. Let
$$w_i=pb_\beta p-2\del p+p(\ta-\tb+4\del\b1+\lam e_\ga)^{1\over 2}(r_i+{1/2}\,\b1)(\ta-\tb+4\del\b1+\lam e_\ga)^{1\over 2}p.$$
Then $w_i\in p\asa p$ and $w_i\to z_{\al\beta}$. Fix an $i$ such that $w_i-z_{\al\beta}\in 1/3\, V$.
Since
$$k-3\del p\leq pb_\beta p-2\del p\leq w_i\leq pa_\al p+2\del p\leq h+3\del p,$$ there is $x$ in $\Cal S$ such that $\|x-w_i\|\leq f(3\del)$.
\enddemo

\proclaim{Corollary 2.3} If $p$ is a closed projection in $A^{**}$ and $h$ is strongly lsc on $p$, then there is $x$ in $p\asa p$ such that $x\leq h$. Moreover, if $t\in(-\infty, 0]$ and $tp\leq h$, then $x$ may be chosen such that $tp\leq x$.
\endproclaim 

\proclaim{Theorem 2.4}Let $p$ be a closed projection in $A^{**}$. If $\fp$ is weak$^*$ metrizable, in particular if $A$ is separable, then $((p\asa p)^{\m})^- =(p\asa p)^{\m}=(p\asa p)^{\sigma}$. Also $((pA_+ p)^{\m})^{-}=(pA_+ p)^{\m} =(pA_+ p)^{\sigma}$.
\endproclaim 

\demo{Proof}This is deduced from Theorem 2.2 in the same way as [B1, Theorem 3.24] and [B1, Corollary 3.25] were deduced from [B1, Theorem 3.15].
\enddemo

\proclaim{Lemma 2.5}If $p$ is a projection in a von Neumann algebra $M$, if $k\leq h$ and $pkp\leq y\leq php$ for $h, k, y\in \msa$, and if $\ep>0$, then there is $x$ in $\msa$ such that $pxp=y$ and $k-\ep\b1\leq x\leq h+\ep\b1$.
\endproclaim 

\demo{Proof}We take $x=k-\ep\b1+z^*z$, where $z=t(h-k+2\ep\b1)^{1\over 2}$ and $\|t\|\leq 1$. It is enough to choose $t$ so that $t(h-k+2\ep\b1)^{1\over 2}p=(y-pkp+\ep p)^{1\over 2}$. Let $q$ be the range projection of $(h-k+2\ep\b1)^{1\over 2}p$. Note that $(h-k+2\ep\b1)^{1\over 2}p$ has closed range, $(h-k+2\ep\b1)^{1\over 2}(p(h-k+2\ep\b1)p)^{-{1\over 2}}$ is a partial isometry from $p$ to $q$, and 
$$q=(h-k+2\ep\b1)^{1\over 2}(p(h-k+2\ep\b1)p)^{-1}(h-k+2\ep\b1)^{1\over 2},$$ where the inverses are taken in $pMp$. So our condition is equivalent to
$$tq=(y-pkp+\ep p)^{1\over 2}(p(h-k+2\ep\b1)p)^{-1}(h-k+2\ep\b1)^{1\over 2},$$
and it suffices to check that with this definition, $\|tq\|\leq 1$. This is done by verifying that $(tq)^*(tq)\leq q$.
\enddemo

\noindent
We state a result from \cite{B1} for convenience.

\proclaim{Lemma 2.6} ([B1, Lemma 3.1 (a)]) If $q$ is a projection in a von Neumann algebra $M$, $s_1\in M$, $\|s_1q\|\leq 1$, and $\|s_1\|\leq 1+\ep$, then there is $s$ in $M$ such that $sq=s_1q,  \|s\|\leq 1$, and $\|s-s_1\|\leq\sqrt{2\ep+\ep^2}$.
\endproclaim 

\proclaim{Lemma 2.7} Assume $p$ is a projection in a von Neumann algebra $M$, $y, h, k\in \msa$, $h\geq k$, $php\geq pyp\geq pkp$, $k\leq y\leq h+\ep\b1$, $p(h-k)p\geq\eta p$, and $0<\ep<\eta \leq \|h-k\|$. Then there is $x$ in $\msa$ such that $pxp=pyp$, $k\leq x\leq h$, and $\|x-y\|\leq C_0({\ep\over \eta})^{1\over 4}\|h-k\|$, where $C_0$ is a universal constant.
\endproclaim 

\demo{Proof}Since
$$0\leq y-k\leq h-k+\ep\b1={(h-k)^{1\over 2}\choose\ep^{1\over 2}\b1}^*{(h-k)^{1\over 2}\choose\ep^{1\over 2}\b1},$$
there are $t_1, t_2$ in $M$ with $\|(t_1\; t_2)\|\leq 1$, such that
$$(y-k)^{1\over 2}=(t_1\; t_2){(h-k)^{1\over 2}\choose\ep^{1\over 2}\b1}=t_1(h-k)^{1\over 2}+\ep^{1\over 2}t_2.$$
We set $x=k+z^*z$, where $z=s(h-k)^{1\over 2}$ and $\|s\|\leq 1$. To achieve $pxp=pyp$, we choose $s$ so that $s(h-k)^{1\over 2}p=(y-k)^{1\over 2}p$. Equivalently, 
$$s(h-k)^{1\over 2}p=t_1(h-k)^{1\over 2}p+\ep^{1\over 2}t_2p, \qquad\text{or}$$
$$(s-t_1)(h-k)^{1\over 2}p=\ep^{1\over 2}t_2p.\leqno(3)$$
Then, by using formulas similar to those relating to $q$ in the proof of 2.5, we see that (3) is equivalent to
$$(s-t_1)q=\ep^{1\over 2}t_2(p(h-k)p)^{-1}(h-k)^{1\over 2},$$
where $q$ is the range projection of $(h-k)^{1\over 2}p$ and the inverse is taken in $pMp$. Also, (3) is equivalent to
$$sq=t_1q+\ep^{1\over 2}t_2(p(h-k)p)^{-1}(h-k)^{1\over 2}=(y-k)^{1\over 2}(p(h-k)p)^{-1}(h-k)^{1\over 2}.\leqno(4)$$
By verifying that $(sq)^*(sq)\leq q$, we see that (4) implies $\|sq\|\leq 1$.

We start with $s_1=t_1+\ep^{1\over 2}t_2(p(h-k)p)^{-1}(h-k)^{1\over 2}$, and note that (4) is satisfied with $s_1$ in place of $s$.
Also, since $(p(h-k)p)^{-{1\over 2}}(h-k)^{1\over 2}$ is  partial isometry, we have
$$\|s_1-t_1\|\leq({\ep\over \eta})^{1\over 2}\quad\text{and}\quad \|s_1\|\leq 1+({\ep\over \eta})^{1\over 2}.\leqno(5)$$
Then by Lemma 2.6, we can find $s$ so that $\|s\|\leq 1$, $sq=s_1q$, and
$$\|s-s_1\|\leq \bigg(2({\ep\over \eta})^{1\over 2}+{\ep\over \eta}\bigg)^{1\over 2}\leq C_1({\ep\over \eta})^{1\over 4}.$$
Thus,
$$\|z-(y-k)^{1\over 2}\|\leq C_1({\ep\over \eta})^{1\over 4}\|h-k\|^{1\over 2}+({\ep\over \eta})^{1\over 2}\|h-k\|^{1\over 2}+\ep^{1\over 2}\leq C_2({\ep\over \eta})^{1\over 4}\|h-k\|^{1\over 2},$$
where in the middle term we have used the estimate (5) for $\|s_1-t_1\|$.\newline
Finally,
$$\|x-y\|\leq\|z-(y-k)^{1\over 2}\|(\|z\|+\|y-k\|^{1\over 2})\leq C_2({\ep\over \eta})^{1\over 4}\|h-k\|^{1\over 2}(2\|h-k\|^{1\over 2}+\ep^{1\over 2})\leq C_0({\ep\over \eta})^{1\over 4}\|h-k\|.$$
\enddemo

\proclaim{Lemma 2.8} Assume $p$ is a projection in a von Neumann algebra $M$, $y,h, k\in \msa$, $h\geq k$, $php\geq pyp\geq pkp$, $k-\ep\b1\leq y\leq h+\ep\b1$, $p(h-k)p\geq \eta p$, and $0<\ep<\eta \leq \|h-k\|$. Then there is $x$ in $\msa$ such that $pxp=pyp$, $k\leq x\leq h$, and $\|x-y\|\leq C({\ep\over \eta})^{1\over 4}\|h-k\|$, where $C$ is a universal constant.  
\endproclaim 

\demo{Proof} First apply Lemma 2.7 with $k-\ep\b1$ in place of $k$. We obtain $x_1$ such that $px_1p=pyp$, $k-\ep\b1\leq x_1\leq h$ and $\|x_1-y\|\leq C_0({\ep\over \eta})^{1\over 4}(\|h-k\|+\ep)$. Then apply the symmetric version of Lemma 2.7 (i.e., apply Lemma 2.7 to $-h\leq -x_1\leq -k+\ep \b1$) to obtain $x$ such that $pxp=pyp$, $k\leq x\leq h$, and $\|x-x_1\|\leq C_0({\ep\over \eta})^{1\over 4}\|h-k\|.$
\enddemo

\example{Remark 2.9} Lemma 3.1 (b) of \cite{B1} asserts the following:

(6)\quad
If $p$ and $q$ are projections in a von Neumann algebra $M$, $t\in M$, $\|ptq\|\leq 1$, and $\|t\|\leq 1+\ep$, then there is $t'$ in $M$ such that $pt'q=ptq$, $\|t'\|\leq 1$ and $\|t'-t\|\leq 2\sqrt{2\ep+\ep^2}$. 
\endexample

\noindent
S. Wassermann [W, p.68] pointed out that the proof of the estimate $2\sqrt{2\ep+\ep^2}$ in \cite{B1} is wrong and that the best estimate proved by the argument in \cite{B1} is $O(\ep^{1\over 4})$ rather then $O(\ep^{1\over 2})$. Lemma 2.8, in the special case $h=\b1$, $k=-\b1$, also gives a $O(\ep^{1\over 4})$ estimate for the special case of (6) where $p=q$ and $t^*=t$. However, (6) is actually correct as stated. Here is a sketch of the proof: Represent $t$ by a matrix $\pmatrix a & b\\
c & d\endpmatrix$ with $a=ptq$, $b=pt(\b1-q)$, etc. The proof of [B1, Lemma 3.1 (a)] produces explicit choices for $b_1$, $c_1$ such that $\|b_1-b\|\leq\sqrt{2\ep+\ep^2}$, $\|c_1-c\|\leq\sqrt{2\ep+\ep^2}$, $\|(a \medspace b_1)\|\leq 1$, and $\|{a\choose c_1}\|\leq 1$. Let $t_1= \pmatrix a & b_1\\
c_1 & d\endpmatrix$. Then $\|t_1-t\|\leq\sqrt{2\ep+\ep^2}$, and it can be shown, by verifying that $t_1^*t_1\leq (1+\ep)^2 \pmatrix q & 0\\
0 & \b1-q\endpmatrix$, that $\|t_1\|\leq 1+\ep$. Then an application of [B1, Lemma 3.1 (a)] to $t_1$ yields $t'=\pmatrix a & *\\
c_1 & *\endpmatrix$ such that $\|t'-t_1\|\leq\sqrt{2\ep+\ep^2}$ and $\|t'\|\leq 1$.

In \cite{B1}, instead of $t_1$ we used  $t_0=\pmatrix a& b\\
c_1 & d\endpmatrix$ and provided a fallaceous proof that $\|t_0\|\leq 1+\ep$. Presumably, the best that can be said is that $\|t_0\|\leq 1+\ep+\sqrt{2\ep+\ep^2}=1+\del$; and thus the argument in \cite{B1} produces a $t'$ with $\|t'-t_0\|\leq \sqrt{2\del+\del^2}$.

The estimate in (6) is sharp to within a factor of 2, but we know nothing about the sharpness of the estimate in Lemma 2.8. However, if we assume in 2.8 that $h-k\geq \eta \b1$, the estimate can be improved to $2\|h-k\|\bigg({\ep\over \eta}+{\ep^2\over \eta^2}\bigg)^{1\over 2}$, even without the hypothesis $\ep<\eta$.

The interpolation--extension theorem will be stated in a very general form. In all the known applications, $q=\b1$.

\proclaim{Theorem 2.10} Assume $p$ and $q$ are closed projections in $A^{**}, p\leq q, h\in((q\asa q)^{\m})^-$, $k\in ((q\asa q)_{\m})^-, h\geq k, p(h-k)p\geq \eta p$ for $\eta >0$, $y\in p\asa p$ and $pkp\leq y\leq pkp$.

(i) There is $x$ in $q\asa q$ such that $pxp=y$ and $k\leq x\leq h$.

(ii) If $x'\in q\asa q, px'p=y$, and $k-\ep q\leq x'\leq h+\ep q$ for $\ep<\eta \leq \|h-k\|$, then $x$ in part (i) can be chosen so that $\|x-x'\|\leq C'({\ep\over \eta})^{1\over 4}\|h-k\|$, where $C'$ is a universal constant.

(iii) If $\Cal A=\{x\in q\asa q: pxp=y\quad\text{and}\quad k\leq x\leq h\}$ and $\cB=\{x\in q\asa^{**} q: pxp=y\quad\text{and}\quad k\leq x\leq h\}$ then $\Cal A$ is $\sigma$--weakly (or equivalently, $\sigma$--strongly) dense in $\Cal B$. 
\endproclaim 

\demo{Proof}We first show
$$\forall \ep>0, \exists x\in q\asa q\; \text{such that }\; pxp=y\; \text{ and }\; k-\ep q\leq x\leq h+\ep q.\leqno(7)$$
To prove (7), note that by Lemmas 2.5 and 2.8 there is $\tx$ in $q\asa^{**} q$ such that $p\tx p=y$ and $k\leq\tx\leq h$.
Then by 2.2 (iii) there is a net $(x_\al)$ in $q\asa q$ such that $k\leq x_\al\leq h$ and $x_\al\to\tx$ $\sigma$--weakly. Since $pA^{**}p$ is the bidual of $pAp$, this implies that $px_\al p\to y$ in the weak Banach space topology of $pAp$. Therefore there is $x_0$, a suitable convex combination of the $x_\al$'s, such that $\|px_0p-y\|\leq\ep$. By 1.2 we can find $z$ in $q\asa q$ such that $\|z\|\leq \ep$ and $pzp=y-px_0p$. Then take $x=x_0+z$.

Next we show

(8)\qquad If $x'\in q\asa q$, $px'p=y, k-\ep q\leq x'\leq h+\ep q$, $0 < \ep<\eta \leq \|h-k\|$, and $\del >0$, then there is $x\in q\asa q$ such that $pxp=y, k-\del q\leq x\leq h+\del q$, and $\|x-x'\|\leq C({\ep\over \eta})^{1\over 4}\|h-k\|$, where $C$ is as in Lemma 2.8.

To prove (8) let $\Cal C=\{x\in q\asa q: k-\del q\leq x\leq h+\del q\}$ and $\cD=\{x\in q\asa q: pxp=y$ and $\|x-x'\|\leq C({\ep\over \eta})^{1\over 4}\|h-k\|\}$. Since $\exists x\in q\asa q$ such that $k\leq x\leq h$, $\Cal C$ has non--empty interior in $q\asa q$. Thus if $\Cal C \cap\cD=\phi$, the Hahn--Banach separation theorem implies the existence of a non--trivial bounded linear functional $f$ on the real Banach space $q\asa q$ such that $\sup\limits_{x\in \Cal C} f(x)\leq\inf\limits_{x\in\cD}f(x)$. 
By Lemma 2.8 there is $\tx$ in $q\asa^{**}q$ such that $k\leq\tx\leq h, p\tx p=y$, and $\|\tx-x'\|\leq C({\ep\over \eta})^{1\over 4}\|h-k\|$. 
Since $\tx$ is in the $\sigma$--weak closure of $\{x\in q\asa q: h\leq x\leq k\}$, we see that $f(\tx)<\sup\limits_{x\in \Cal C} f(x)$. Let $L=\{x\in A:xp=0\}$ and $L^*=\{x\in A:px=0\}$. 
Then $\cD=x'+qBq$, where $B$ is a ball in $(L+L^*)_{\sa}$. 
Since $L+L^*$ is closed by \cite{C1}, the $\sigma$--weak closure of $(L+L^*)$ can be identified with its bidual; and we conclude that $\tx$ is in the $\sigma$--weak closure of $\cD$. therefore $f(\tx)\geq \inf\limits_{x\in \cD}f(x)$, a contradiction.

Now it is routine to combine (7) and (8) to prove (i). And (ii) follows from (8), where $C'$ can be any number greater than $C$. For (iii), given $x$ in $\cB$, there is a net $(x_\al)$ in $q\asa q$ such that $k\leq x_\al\leq h$ and $x_\al\to x$, $\sigma$--weakly, as above. And as above, $px_\al p\to y$ in the weak Banach space topology of $pAp$. Thus there is a net $(x'_\beta)$, consisting of convex combinations of the $x_\al$'s such that $x'_\beta\to x$ $\sigma$--weakly and $px'_\beta p\to y$ in norm. Then we can find $z_\beta$ in $q\asa q$ such that $z_\beta\to 0$ in norm and $pz_\beta p=y-px'_\beta p$. So if $x''_\beta=x'_\beta+z_\beta$, we have $x''_\beta\to x$ $\sigma$--weakly, $px''_\beta p=y$, and $k-\del_\beta q\leq x''_\beta\leq h+\del_\beta q$, where $\del_\beta\to 0$. Then by applying part (ii) to $x''_\beta$, we find $x'''_\beta$ in $\Cal A$ such that $\|x'''_\beta-x''_\beta\|\to 0$. 

If $p$ is central, the hypothesis $p(h-k)p\geq \eta p$ can be omitted from theorem 2.10. One can show in this case that $\th=h(q-p)+y$ and $\tk=k(q-p)+y$ give an lsc and a usc function on $F(q)$. Then apply Theorem 2.2 to $q, \th$ and $\tk$. But in general, this hypothesis cannot be omitted, even if $A$ is unital $q=\b1$, and $p,h, k$ are all in $A$.
\enddemo

\example{Example 2.11} Let $A=c\otimes\bM_2$, the algebra of convergent sequences in $\bM_2$. Let $(t_n)$ be a non--convergent sequence in (0,1), and let $(\del_n)$ be a sequence such that $\del_n>0$ and $\del_n\to 0$. Let $p$ be the constant sequence $\pmatrix 1& 0\\
0 & 0\endpmatrix$, and let $h$ and $k$ in $\asa$ be given by $h_n=\pmatrix t_n\del_n & \del_n^{1\over 2}\\
\del_n^{1\over 2} & {1\over 2}\endpmatrix$, $k_n=\pmatrix (t_n-1)\del_n & 0\\
0 & -{1\over 2}\endpmatrix$. If $y=0$, in $p\asa p$, then all hypotheses of 2.10 except $p(h-k)p\geq \eta p$ are satisfied, and also $(\b1-p)(h-k)(\b1-p)\geq \b1-p$. But there is no $x$ in $\asa$ such that $pxp=y$ and $k\leq x\leq h$. If $x$ existed, then $x_n=k_n+s_n(h_n-k_n)$, $s_n\in[0,1]$, since $h_n-k_n$ has rank 1. Since $p_nx_np_n=0$, we see that $s_n=1-t_n$. Then $x_n=\pmatrix 0 & *\\
* & {1\over 2}-t_n\endpmatrix$, which is absurd, since $({1\over 2}-t_n)$ is not convergent.
\endexample                                                

\proclaim{Proposition 2.12} Assume $p$ is a closed projection in $A^{**}$, $h\in((p\asa p)^{\m})^-$, and $sp\leq h\leq tp$ with $s\leq 0$, $t\geq 0$. 
If $\fp$ is weak$^*$ metrizable, in particular if $A$ is separable, then there is $\th$ in $((\asa)^{\m})^-$ such that $p\th p=h$ and $s\b1\leq\th\leq t\b1$. Under the same hypothesis, $p((\Asa)^{\m})p=(p\Asa p)^{\m}$.
\endproclaim 

\demo{Proof} Use 2.3 to find $x$ in $p\asa p$ such that $sp\leq x\leq h$. Then, since by 2.4 $h-x$ is in $(pA_+p)^\sigma$, there is an increasing sequence $(x_n)$ in $p\asa p$ such that $x_n\to h$ and $x_n\geq sp$. Let $y_n=(1-{1\over n})x_n+{1\over n}sp$. Then $(y_n)$ is increasing, and, except in the trivial case $s=t=0, tp-y_n\geq \del_n p$ for $\del_n>0$. Thus we can use Theorem 2.10 to construct recursively $a_1, a_2,\dots$ in $\asa$ such that $pa_n p=(1-{1\over n})x_n$ and 
$$s\b1\leq a_n+{s\over n}\b1\leq a_{n+1}+{s\over n+1}\b1\leq t\b1.$$
Then, since ${s\over n}\to 0$, all conditions are met if we take $\th=\lim(a_n+{s\over n}\b1)$.

The last sentence follows easily.
\enddemo

\example{Remark 2.13} Although for each $h$ in $(p\Asa p)^{\m}$ there is $\th$ in $(\Asa)^{\m}$ such that $p\th p=h$, $\|\th\|$ may have to be much larger than $\|h\|$. (It doesn't help to relax the requirement to $\th\in((\Asa)^{\m})^-$.) In particular, if $h$ is in $((p\asa p)^{\m})^-$ and $sp\leq h\leq tp$ for $s>0$, there need not exist $\th$ in $((\asa p)^{\m})^-$ such that $p\th p=h$ and $(s-\ep)\b1\leq\th\leq(t+\ep)\b1$. However, there does exist $\th$ with $(s-\ep)\b1\leq\th$.

\endexample
 
\S 3. Continuous elements

Each type of semicontinuity gives rise to a concept of continuity, where $h$ is continuous if it is both lsc and usc. 
We already know that $h$ is strongly continuous on $p$, in this sense, if and only if $h\in p\asa p$ if and only if $h$ is in $A_0(p)$. 
For the other two types of semicontinuity, when $p=\b1$, it was proved in \cite{AP} that $(\Asa)^{\m}\cap(\Asa)_{\m}=M(A)_{\sa}$ and $((\Asa)^{\m})^-\cap((\Asa)_{\m})^-=QM(A)_{\sa}$, where $QM(A)$ is the quasi--multiplier  space of $A$. 
Since in general $pAp$ is not an algebra, we cannot characterize the middle and weakly continuous elements of $p\asa^{**}p$ in terms of multiplier properties. 
Instead we will show that they are related to $pM(A)_{\sa}p$ and $pQM(A)_{\sa}p$.

\noindent
Most of the content of the next lemma is needed to deal with the non--separable case.
 
\proclaim{Lemma 3.1} Assume $p$ is a closed projection in $A^{**}$ and $A$ is $\sigma$--unital. If $h\in(p\Asa p)^\sigma$, then for sufficiently large $\lam$, there is an increasing sequence $(x_n)$ in $pA_+p$ such that $x_n\to h+\lam p$ and $h+\lam p-x_n\geq\del_np$ with $\del_n>0$.
\endproclaim  

\demo{Proof} Assume $h$ is the limit of an increasing sequence $(\lam_np+y_n)$, $\lam_n\in\bR$, $y_n\in p\asa p$. By Lemma 1.4, $\{\lam_n\}$ is bounded if $p$ is not compact; and if $p$ is compact, we may also arrange that $\{\lam_n\}$ is bounded. Choose $\lam$ large enough that $\mu_n=\lam+\lam_n>0$, $\forall n$, and that $\mu_np+y_n\geq p$. Thus $\mu_n\in (0,M]$ for some $M>0$. Let $\cP=\{a\in A_+:\|a\|<1\}$, as in the proof of 1.3, and let $(e_n)$ be a sequential approximate identity for $A$ with $e_n\in\cP$. Then, using the same technique as in the proof of 1.3, we can recursively construct $f_1,f_2,\dots$ in $p\cP p$ such that
$$0 \leq \mu_nf_n+y_n\leq  \mu_{n+1}f_{n+1}+y_{n+1},\quad\text{and}\quad f_n\geq pe_n p.$$
We claim that the choice $x_n=\mu_nf_n+y_n$ needs all requirements. Since $(x_n)$ is bounded and increasing, it converges to some $h'$ in $p\asa^{**}p$, and clearly $h'\leq h+\lam p$. But since $x_n\geq(\lam+\lam_n)p+y_n-M(p-pe_n p)$ and $pe_np\to p$, we also have $h'\geq h+\lam p$.
\enddemo

\proclaim{Theorem 3.2} Assume $A$ is a $\sigma$--unital $C^*$--algebra, $p$ is a closed projection in $A^{**}$, and $h\in p\asa^{**}p$.

(i) Then $h\in pM(A)_{\sa} p$ if and only if $h\in (p\Asa p)^\sigma\cap (p\Asa p)_\sigma$.

(ii) If $\fp$ is weak$^*$ metrizable, in particular if $A$ is separable, then $pM(A)_{\sa} p=(p\Asa p)^{\m}\cap(p\Asa p)_{\m}$.

(iii) If $\fp$ is weak$^*$ metrizable, in particular if $A$ is separable, then $(pQM(A)_{\sa} p)^-=((p\Asa p)^{\m})^{-}\cap((p\Asa p)_{\m})^-$. 
\endproclaim  
 
\demo{Proof} (i) Since $A$ is $\sigma$--unital, $M(A)_{\sa}=(\Asa)^\sigma\cap(\Asa)_\sigma$. 
Thus the necessity of the condition is clear. 
So assume $h\in (p\Asa p)^\sigma\cap(p\Asa p)_\sigma$. 
Choose $\lam$, $\mu>0$ so that $h+\lam p$ is the limit of an increasing sequence $(x_n)$ with the properties in 3.1 and $h-\mu p$ is the limit of a decreasing sequence $(y_n)$ with symmetrical properties. 
Let $e$ be a strictly positive element of $A$. 
We will recursively construct sequences $(a_k)$ and $(b_k)$ in $A_{\sa}$ such that:

(a) $-\lam\b1+a_k\leq -\lam\b1+a_{k+1}\leq \mu\b1+b_{k+1}\leq \mu\b1+b_k$,

(b) $pa_kp\in {\text{co}}\{x_n:n\geq k\}$ and $pb_k p\in {\text{co}}\{y_n:n\geq k\}$, where ${\text{co}}$ denotes convex hull, and 

(c) $\|e(\mu\b1+b_k+\lam\b1-a_k)e\|\leq{1\over k}, k\geq 2$.

Choose $a_1$ arbitrarily such that $pa_1p=x_1$. Then use Theorem 2.10 to construct $b_1$ such that $pb_1p=y_1$ and $b_1\geq(-\lam-\mu)\b1+a_1$.

Now assume $a_1,\dots,a_k$ and $b_1,\dots,b_k$ have already been constructed. Choose $N\geq k+1$ so that $pa_kp\in {\text{co}}\{x_n:n<N\}$ and $pb_kp\in {\text{co}}\{y_n:n<N\}$. We first construct $a'_n, b'_n$ in $\asa$, for $n\geq N$ so that:
$$(a')\quad -\lam\b1+a_k\leq -\lam\b1+a'_N\leq -\lam\b1+a'_{N+1}\leq\cdots\leq \mu\b1+b'_{N+1}\leq \mu\b1+b'_N\leq \mu\b1+b_k,\; \text{and}$$

\noindent
$(b')\quad pa'_np=x_n, pb'_np=y_n.$

\noindent
This is just a matter of successive applications of Theorem 2.10. (Choose first $a'_N$, then $b'_N$, then $a'_{N+1}$, etc.) Then let $f_n=(\lam+\mu)\b1+b'_n-a'_n$ and $f=\lim f_n$. Then $pfp=0$. Since $0\leq f\leq f_n$, Theorem 2.10 (iii) produces a net $(z^{(n)}_\al)_{\al\in D_n}$ such that $z^{(n)}_\al\in\asa$, $0\leq z^{(n)}_\al\leq f_n$, $pz^{(n)}_\al p=0$ and $z^{(n)}_\al\to f$ $\sigma$--weakly. Let ${\Cal S}=\{e(f_n-z^{(n)}_\al)e: n\geq N, \al\in D_n\}$. Since $f_n\to f$ $\sigma$--weakly and $z^{(n)}_\al\to f$ $\sigma$--weakly for each $n$, we see that 0 is in the $\sigma$--weak closure of ${\Cal S}$. But ${\Cal S}\subset A$ and the restriction to $A$ of the $\sigma$--weak topology is the weak Banach space topology of $A$. Thus there is $v$ in ${\text{co}}({\Cal S})$ such that $\|v\|<{1\over k+1}$. Then there are $n_1,\dots,n_l\geq N$ and $s_1,\dots, s_l\geq 0$ such that $\sum s_i=1\;\text{and}\; v=e(f'-z)e$ where $f'=\sum s_if_{n_i}, z\in A_+, pzp=0$, and $z\leq f'$, Then we can take $a_{k+1}=z+\sum s_ia'_{n_i}$ and $b_{k+1}=\sum s_ib'_{n_i}$.

Now it is clear that $\lim(-\lam\b1+a_k)=\lim(\mu\b1+b_k)$. if $\th$ is this limit, then $\th\in M(A)_{\sa}$ by \cite{AP} and $p\th p=h$. 

(ii) Follows from (i) and Theorem 2.4.

(iii) Since $QM(A)_{\sa}=((\Asa)^{\m})^-\cap((\Asa)_{\m})^-$ by \cite{AP}, it is clear that $pQM(A)_{\sa}p\subset ((p\Asa p)^{\m})^-\cap((p\Asa p)_{\m})^-$.
\noindent
Thus also $(pQM(A)_{\sa} p)^-\subset ((p\Asa p)^{\m})^-\cap((p\Asa p)_{\m})^-$. 
For the converse, take $f\in((p\Asa p)^{\m})^-\cap((p\Asa p)_{\m})^-$ and $\ep>0$. 
Choose $h\in(p\Asa p)^{\m}$ such that $\|h-f\|<\ep$, and let $h_1=h+\ep p$. 
Then $h_1\geq f$, $\|h_1-f\|<2\ep$, and $h_1-f\geq\del p$ for some $\del>0$. 
In a symmetrical way, find $k_1$ in $(p\Asa p)_\m$ with $k_1\leq f$. 
By 2.12 there are $\th$ in $(\Asa )^\m)$ and $\tk$ in $(\Asa)_\m$ such that $p\th p=h_1$ and $p\tk p=k_1$. 
Since $p(\th-\tk)p\geq\del p$, for $\lam$ sufficiently large we have $\th+\lam(1-p)\geq\tk$. 
Then by [B1, Theorem 3.26 (c)] there is $\tf$ in $QM(A)_{\sa}$ such that $\tk\leq\tf\leq\th+\lam(1-p)$. Then $\|p\tf p-f\|<2\ep$.
\enddemo 

\example{Remark 3.3} (i) The Non--commutative Tietze Extension Theorem was proved in the separable case in \cite{APT} and in general form in [P3, Theorem 10]. It states that the natural map from $M(A)$ to $M(A/I)$ is surjective when $A$ is $\sigma$--unital and $I$ is a closed two--sided ideal. Part (i) of the theorem specializes to this when $p$ is in the center of $A^{**}$. Another theorem that specializes to [P3, Theorem 10] is [B1, Theorem 3.43(b)]. It implies, when $A$ is $\sigma$--unital and $h$ is in $p\asa^{**}p$ that $h=p\th$ for some $\th$ in $M(A)$ such that $p\th=\th p$ if and only if $h$ is $q$--continuous on $p$. It is interesting to note that part (i) of Theorem 3.2 was proved with Pedersen's techniques, whereas [B1, Theorem 3.43] was proved by totally different methods (but still related to semicontinuity).

(ii) It was shown in [B1, Example 3.13] that $pM(A)p$ need not be norm closed. 
Therefore when we wish to find $\th$ in $M(A)_{\sa}$ such that $p\th p=h$, we may need to take $\|\th\|$ much larger than $\|h\|$.

(iii) Assume $\fp$is metrizable. There are two natural vector spaces intermediate between $pM(A)_{\sa}p$ and $(pQM(A)_{\sa}p)^-$, namely $(pM(A)_{\sa}p)^-$ and $pQM(A)_{\sa} p$. In general these spaces are not comparable. In fact [B1, Example 3.13] showed that $(pM(A)_{\sa}p)^-$ need not be contained in $pQM(A)_{\sa}p$. And it is easy to find examples, with $p=\b1$, where $QM(A)_{\sa}\not\subset M(A)_{\sa}=(M(A)_{\sa})^-$. By taking the direct sum of two examples, we get non--comparability.

In [B1, Proposition 2.3] we showed that $(\Asa)^{\m}\cap((\Asa)_{\m})^-= ((\Asa)^{\m})^-\cap(\Asa)_{\m}=M(A)_{\sa}$, but we do not know whether $(p\Asa p)^{\m}\cap((p\Asa p)_{\m})^-=((p\Asa p)^{\m})^-$\break
$\cap(p\Asa p)_\m$. 
These sets are also intermediate between $pM(A)_{\sa}p$ and $(pQM(A)_{\sa}p)^-$. If they are equal, they are both equal to $pM(A)_{\sa}p$; but if they are unequal, neither is a vector space, since they are negatives of one another.

(iv) Theorem 2.2 is a satisfactory generalization to closed faces of [B1, Theorem 3.15], the ``strong interpolation theorem''. There are less satisfactory versions for closed faces of [B1, Theorem 3.26 (c)] and [B1, Theorem 3.40], the ``middle'' and ``weak'' interpolation theorems. These are just corollaries of the results in \cite{B1}. 

The middle version states that if $A$ is $\sigma$--unital, $h$ is $q$--lsc on $p$, $k$ is $q$--usc on $p$, and $\overset q\to{h \geq k}$, then there is $f$ in $pM(a)_{\sa} p$ such that $k\leq f\leq h$. This is proved by applying [B1, Theorem 3.40] to $h+t(\b1-p)$ and $k+s(\b1-p)$, for $t$ sufficiently large and $s$ sufficiently small, to obtain $\tf$ in $M(A)_{\sa}$ such that $k+s(\b1-p)\leq\tf\leq h+t(\b1-p)$. 
Then take $f=p\tf p$. 
Although this result is in some sense an exact analogue of [B1, Theorem 3.40], we would really like to have $f$ $q$--continuous on $p$. 
This follows from the stated conclusion when $p=\b1$ but not in general.

The weak version states that if $A$ is $\sigma$--unital, $\fp$ is weak$^*$ metrizable $h$ is weakly lsc on $\fp$, $k$ is weakly usc on $\fp$, and $h-k\geq\del p$ for $\del>0$, then there is $f$ in $pQM(A)_{\sa}p$ such that $k\leq f\leq h$. This is deduced from [B1, Theorem 3.26(c)] similarly to the proof of part(iii) of Theorem 3.2. (Approximate $h$ from below by $(p\Asa p)^{\m}$ and $k$ from above by $(p\Asa p)_{\m}$.)

\endexample
 
\example{Example 3.4} We are now able to provide the example promised in Remark 1.6, an element of $p\asa^{**}p$ which is an lsc functional on $\fp\cap S$ but is not weakly lsc on $\fp$. Let $A=c\otimes \bK$ and let $p$ be as in [B1, Example 3.13]. Thus let $v_{k,n}=(1-1/k)^{1\over 2}e_k+(1/k)^{1\over 2} e_{n+k}, k=1,\dots,n$, where  $e_1,e_2,\dots$ is the standard basis of $l^2$. Then $p$ is given by $(p_n)$ where $p_n=\sum^n_{k=1}\;v_{k,n}\times v_{k,n}$ for $n<\infty$ and $p_\infty=\b1$. The fact that $p_\infty=\b1$ implies $p$ is closed.

We claim that the restriction to $\fp\cap S$ of the weak$^*$ topology agrees with the norm topology, which implies that every element of $B_0(p)$ is continuous on $\fp\cap S$. Note that each $\var$ in $Q(A)$ is given by $\{\var_n: 1\leq n\leq\infty\}$ where $\var_n\in Q(\bK)$ and $\|\var\|=\|\var_\infty\|+\sum^\infty_1\|\var_n\|$. We will represent $\var$ by a pair $(\var',\var'')$, where for $h$ in $A^{**}$, $\var'(h)=\sum^\infty_1\var_n(h_n)$ and $\var''(h)=\var_\infty(h_\infty)$. Now suppose $\var_i$ and $\var$ are in $\fp\cap S$ and $\var_i\to\var$ weak$^*$. It is permissible to pass to a subsequence so that $\var'_i\to\theta$ and $\var''_i\to\psi$ for some $\theta,\psi$. Clearly $\psi=(0,\psi'')$, and we proceed to show that $\theta=(\theta', 0)$. For temporarily fixed $k,l$ define $h$ in $QM(A)$ by
$$h_n=\cases 0, &n<\max(k,l),\\
e_k\times e_l, & n=\infty,\\
e_k\times e_l-(1-1/k)^{1\over 2}(1-1/l)^{1\over 2}k^{1\over 2}l^{1\over 2} e_{n+k}\times e_{n+l}, & \text{otherwise,}\endcases$$
Since $h$ is weak$^*$ continuous on $S$ and $\|\var'_i\|\to\|\theta\|$, $\theta(h)=\lim\var'_i(h)$. But $\var'_i(h)=0$, since $p_nh_np_n=0$ for $n<\infty$; and $\theta(h)=\theta''(h)$ for the same reason. Therefore $\theta''=0$, whence $\theta=\var'$ and $\psi =\var''$. Since we now have that $\var_{in}\to\var_n$ for $1\leq n\leq\infty$ and $\sum^\infty_{n=1}\|\var_{in}\|\to \sum^\infty_{n=1}\|\var_n\|$, and since it is well known that the restriction of the weak$^*$ topology of $\bK^*$ to $S(\bK)$ agrees with the norm topology, it follows that $\|\var_i-\var\|\to 0$, and the claim is proved.

Next we produce an $h$ in $p\asa^{**}p\setminus (pQM(A)_{\sa} p)^-$. It was shown in [B1, Example 3.13(i)] that if $h$ is in $pQM(A)p$ and $h_\infty=0$, then $\lim\sup_{n\to\infty}|(h_nv_{k,n},v_{k,n})|={\text{O}}(1/k^{1\over 2})$. Now it is easily seen that if $h\in(pQM(A) p)^-$ and $h_\infty=0$, then $h$ is in $\{pxp: x\in QM(A)\;\text{and}\; x_\infty=0\}^-$. 
It follows that $\lim_{k\to\infty}\lim\sup_{n\to\infty}|(h_nv_{k,n}v_{k,n})|$\break
$=0$. To find an $h$ not satisfying this condition, let $h_n=-p_n$ for $n<\infty$ and $h_\infty=0$.
 
By Theorem 3.2(iii) $h$ is not both weakly lsc and weakly usc. Since $h$ is obviously usc, it is the promised example of an element of $p\asa^{**}p$ which is an lsc functional on $\fp\cap S$ but is not weakly lsc.
\endexample

\S 4. Semicontinuity and the continuous functional calculus.

In this section we consider the continuous functional calculus within the algebra $pA^{**}p$ for functions which are operator convex, and sometimes also operator monotone. If $f$ is a continuous function on an interval $I$, $f$ is called {\it{operator}} {\it{monotone}}, and $-f$ is called {\it{operator}} {\it{decreasing}}, if for self--adjoint operators $h_1$ and $h_2$ with $\sigma(h_i)\subset I$, $h_1\leq h_2$ implies $f(h_1)\leq f(h_2)$. Also $f$ (still continuous on $I$) is called {\it{operator}} {\it {convex}}, and $-f$ is called {\it {operator}} {\it {concave}}, if for self--adjoint $h_1$ and $h_2$ with $\sigma(h_i)\subset I$ and $t\in[0,1]$, we have $f(th_1+(1-t)h_2)\leq tf(h_1)+(1-t)f(h_2)$, and $f$ is {\it{strongly}} {\it{operator}} {\it{convex}} if the function $f(x_0 + \cdot)$, for $x_0$ in $I$, satisfies the six equivalent conditions of [B1, Theorem 2.36].  One of these conditions is that $pf(php)p \leq f(h)$ whenever $p$ is
 a projection and $h$ is a self-adjoint operator with $\sigma(h) \subset I$, and Ch. Davis proved in [D1] that operator convexity
 is equivalent to the weaker condition $pf(php)p \leq pf(h)p$.  (The translation of the independent variable in $f$ is needed only because one of the six conditions in [B1, 2.36] doesn't make sense unless $0 \in I$.  The symbol $pf(php)p$ can easily be interpreted to make sense, and lead to correct characterizations, even if $0 \notin I$.)
 A self--contained and efficient exposition of operator monotonicity and convexity, from the point of view of operator algebraists, can be found in \cite{HP}. Historical background  can be found in \cite{HP} and \cite{D2}.

If $f$ is operator convex, then for each $x$ in $I$
$$
f(x)=ax^2+bx+c+\int_{t<I}{(x-x_0)^2\over (x-t)(x_0-t)^2}d\mu_-(t)+\int_{t>I}{(x-x_0)^2\over (t-x)(t-x_0)^2}d\mu_+(t).\leqno(9)
$$
Here $a\geq 0$, $x_0$ can be any point in the interior of $I$ (the integrands are obtained by subtracting from $\pm 1/(t-x)$ its first degree Taylor polynomial at $x_0$), and $\mu_\pm$ are positive measures which are finite on bounded sets and such that $\int 1/(1+|t|^3)d\mu_{\pm}(t)<\infty$. 
These conditions are sufficient to imply the convergence of the integrals for $x$ in the interior of $I$. 
If $I$ contains one or both endpoints, then convergence of (9) at such endpoint(s) imposes an additional restriction on $\mu_\pm$. 
Any choice of $a, b,c$, and $\mu_\pm$ meeting all these requirements does indeed produce an operator convex function.

If $f$ is strongly operator convex, then for each $x$ in $I$
$$ f(x) = c + \int_{t<I} {1\over{x-t}}d\mu_{-}(t) + \int_{t>I} {1\over{t-x}}d\mu_{+}(t). \leqno(10)$$
Here $c\geq 0$ and $\mu_\pm$ are positive measures such that $\int 1/(1+|t|)d\mu
_\pm < \infty$.  Again if $I$ contains one or both endpoints, then convergence 
at such endpoint(s) imposes an additional restriction on $\mu_\pm$, and all the conditions imply that (10) does indeed produce a strongly operator convex 
function.  (The measures $\mu_\pm$ appearing for $f$ in (10) are the same ones appearing for $f$ in (9).)

If $f$ is operator monotone, then for each $x$ in $I$
$$f(x)=ax+b+\int_{t<I}{x-x_0\over(x-t)(x_0-t)}d\mu_{-}(t)+\int_{t>I}{x-x_0\over(t-x)(t-x_0)}d\mu_{+}(t).\leqno(11)$$
Here $a\geq 0$, $x_0$ can be any point in $I$ (the integrand is $1/(t-x)-1/(t-x_0))$ and $\mu_\pm$ are positive measures such that $\int 1/(1+t^2)d\mu_{\pm}(t)<\infty$. Again, if $I$ contains one or both endpoints, then convergence at such endpoint(s) imposes an additional restriction on $\mu_\pm $ and all the conditions imply that (11) does indeed produce an operator monotone function.

\proclaim{Lemma 4.1} Let $p$ be a closed projection in $A^{**}$, where $A$ is unital, and let $h$ be in $p\asa p$. 

(i) Then $h^2$ is usc on $p$.

(ii) If $h\geq\ep p$ for some $\ep>0$, then $h^{-1}$ is usc in $A^{**}$ and a fortiori usc on $p$.
\endproclaim 

\demo{Proof} (i) If $h=pap$, then $h^2=papap$. It is routine to deduce from $p\in (\asa)_\m$ that $apa\in (\asa)_\m$. Therefore $p(apa)p\in(p\asa p)_\m$.

(ii) By 1.2 we may write $h=pap$ for $a\geq \ep\b1$. 
Let $y=a^{1\over 2} p$. Since $y^*y\geq\ep p$, $y$ has closed range. 
The range projection, $q$, of $y$ is also the range projection of $yy^*=a^{1\over 2} pa^{1\over 2}$. 
Note that $yy^*$ is usc (as above) and that $\sigma(yy^*)$ omits the interval $(0,\ep)$. 
It is known (cf. [B1, Proposition 2.44(b)] that this implies $q$ is closed. Since (as in the proof of 2.5) $q=a^{1\over 2}(pap)^{-1} a^{1\over 2}$, where the inverse is taken in $pA^{**}p$, and since $a^{1\over 2}$ is invertible, we deduce from the fact that $q$ is usc that also $(pap)^{-1}$ is usc. 
\enddemo

\proclaim{Lemma 4.2} Assume $p$ is a closed projection in $A^{**}$, $A$ is unital, $f$ is an operator convex function on an interval $I$, $h\in p\asa p$, and $\sigma_{pA^{**}p}(h)\subset I$. Then $f(h)$ is usc on $\fp$.
\endproclaim 

\demo{Proof} We use the integral representation of $f$. Then $f(h)$ can be obtained by substituting $h$ for $x$ in (9), thus obtaining a Bochner integral. Since $(p(\asa)_{\m}p)^-$ is a closed cone, it is enough to verify that each term and each value of the integrand is usc on $p$; and this follows from Lemma 4.1.
\enddemo

\proclaim{Theorem 4.3} Assume $p$ is a closed projection in $A^{**}$, $f$ is an operator convex function on an interval $I$, $h\in p\asa p$, and $\sigma_{pA^{**}p}(h)\subset I$. 

(i) If either $p$ is compact or $0\in I$ and $f(0)\leq 0$, then $f(h)$ is strongly usc on $p$.

(ii) If $0\in I$, then $f(h)$ is middle usc on $p$.

(iii) If $0$ is an endpoint of $I$, then $f(h)$ is weakly usc on $p$. 

\noindent
Conversely, if $f$ is a function which satisfies the conclusion of (i) for all closed faces in the unital case, then $f$ is operator convex.
\endproclaim 

\demo{Proof} (i) We apply Lemma 4.2 to $\tA$, identifying $(\tA)^{**}$ with $A^{**}\oplus\bC$. If $p$ is compact, then $p\oplus 0$ is closed and $f(h\oplus 0)$, computed in $(p\oplus 0)(\tA)^{**}(p\oplus 0)$, is the same as $f(h)\oplus 0$. Since the weak$^*$ topologies from $\tA$ and $A$ agree on $\fp$, 4.2 implies $f(h)$ is strongly usc on $p$. If $p$ is not compact, then $\tp=p\oplus 1$ is closed. Applying 4.2 to $h\oplus 0$ (in $\tp\Asa\tp$), we find that $f(h)\oplus f(0)$ is usc on $\tp$, whence $(f(h)-f(0)p)\oplus 0$ is also usc on $\tp$. Since $\fp$ can be identified with $F(\tp)\cap S(\tA)$, this implies $f(h)-f(0)p$ is strongly usc on $p$. And since $f(0)\leq 0$, it follows that $f(h)$ is also strongly usc.

(ii) Follows from (i).

(iii) We may assume $0$ is the left  endpoint of $I$.  There are sequences $(\del_n)$ in $(0,\infty)$ and $(\theta_n)$ in (0,1) such that $\del_n\to 0, \theta_n\to 1$ and $f(\del_n p+\theta_n h)$ is defined for all $n$. Since $f(\del_n p+\theta_n h)$ is obtained by applying $f(\del_n+\cdot)$ to $\theta_n h$, $f(\del_n p+\theta_n h)\in (p\Asa p)_{\m}$ for each $n$. Therefore $f(h)\in ((p\Asa p)_{\m})^-$.

For the converse let $A=c\otimes\bM_\m$, where $\m=k+l$. Matrices will be written in block form, $\pmatrix a &b\\
c &d\endpmatrix$, where $a$ is $k\times k$, $b$ is $k\times l$, etc. Elements of $A^{**}$ will be represented by bounded indexed collections $(h_n)_{1\leq n\leq\infty}$. Let $p$ be the closed projection given by $p_n=\pmatrix \b1_k &0\\
0 &0\endpmatrix$ for $n<\infty$ and $p_\infty=\pmatrix \b1_k &0\\
0 &\b1_l\endpmatrix$ If $h$ is in $p\asa^{**} p$, then for $n$ finite, $h_n$ can be regarded as a $k\times k$ matrix $a_n$, and $h_\infty=\pmatrix a_\infty &b_\infty\\
b_\infty^* &c_\infty\endpmatrix$. It is fairly easy to see that $h$ is in $p\asa p$ if and only $a_n\to a_\infty$ and $h$ is usc if and only if $a'\leq a_\infty$ for every cluster point $a'$ of $(a_n)$. (The last condition is equivalent to: For each $\ep>0$, there is $N$ such that $n\geq N\Rightarrow a_n\leq a_\infty +\ep\b1_{k}$.) Choose a matrix $h_\infty=\pmatrix a_\infty &b_\infty\\
b_\infty^* &c_\infty\endpmatrix$ such that $\sigma(h_\infty)\subset I$ and let $a_n=a_\infty$ for $n<\infty$. Then $h$ is in $p\asa p$, and $f(h)$ is usc if and only if $f(pr(h_\infty))=f(a_\infty)\leq pr(f(h_\infty))$. Here $pr\bigg(\pmatrix a &b\\
c &d\endpmatrix\bigg)=a$. This last condition demanded for all $k,l$, is equivalent to operator convexity by \cite{D1}.
\enddemo
 
\example{Remark 4.4} In [B1, propositions 2.34 and 2.35(b)] another characterization of operator convexity was given; namely, $f$ is operator convex if and only if $f(h)$ is weakly lsc when $h\in QM(A)_{\sa}$, and $\sigma(h)\subset I$ (demanded for all  $C^*$--algebras). It may seem strange that in the one situation operator convexity produces lower semicontinuity and in the other upper semicontinuity. However, if we keep in mind the closure property of $C^*$--algebras under the continuous functional calculus, then a moment's thought convinces us that it all makes sense. None of $QM(A), pAp$, or $pM(A)p$ are $C^*$--algebras in general, but $A$ and $M(A)$ are. (We will treat $pM(A)p$ below.)
\endexample 

We now show that the conditions relating to $0$ and $f(0)$ cannot be dropped from Theorem 4.3. Of course, $f(0)$, computed by the continuous functional calculus for $0$ in $p\asa p$, is $f(0)p$; and $f(0)p$ is strongly usc if and only if $p$ is compact or $f(0)\leq 0$. Also, there does not exist $h$ in $p\asa p$ with $0\not\in\sigma_{pA^{**}p}(h)$ unless $p$ satisfies the condition (1) from \S1.

\example{Example 4.5} Here $A=c\otimes\bK$, as in Example 1.8, and elements of $A^{**}$ are identified with bounded collections $(h_n)_{1\leq n\leq\infty}$, where $h_n\in B(l^2)$. Choose $\theta$ in $(0, {\pi\over 2})$ and let $v_n=\cos\theta e_1+\sin\theta e_{n+2}$. Then define $p_n$ for $n<\infty$ as $v_n\times v_n+e_2\times e_2$ and $p_\infty$ as $e_1\times e_1+e_2\times e_2$. Then $p=(p_n)$ is closed. For $h=(h_n)$ in $p\asa ^{**}p$, we will represent $h_n$ as a $2\times 2$  matrix $\pmatrix a_n &b_n\\
\overline{b_n} &c_n\endpmatrix$ relative to the basis $\{v_n, e_2\}$, for $n<\infty$, and $h_\infty$ as a matrix $\pmatrix a_\infty &b_\infty\\
\overline{b_\infty} &c_\infty\endpmatrix$ relative to the basis $e_1, e_2$. Let $r=\pmatrix s &t\\
\ot &u\endpmatrix$ be a positive matrix such that $tr(r)=s+u\leq 1$, and let $r'=\pmatrix s\cos^2\theta &t\cos\theta\\
\ot\cos\theta &u\endpmatrix$. Let $\var_n$ in $\fp$ be given by $\var_n(h)=tr(r h_n)$ for $n<\infty$, and define $\var_\infty$ by $\var_\infty(h)=tr(r'h_\infty)$. Then $\var_n\to\var_\infty$ weak$^*$. Thus if $h$ is strongly usc, $tr(h_\infty r')\geq\lim\sup tr(h_n r)$. This implies that $\pmatrix a_\infty\cos^2\theta &b_\infty\cos\theta\\
\overline{b_\infty}\cos\theta &c_\infty\endpmatrix\geq \pmatrix a &b\\
\ob &c\endpmatrix$ for every cluster point $\pmatrix a &b\\
\ob &c\endpmatrix$ of $(h_n)$. (This necessary condition is also sufficient for strong upper semicontinuity.) It follows that if $h$ is middle usc, then either $c_\infty>c$ or $c_\infty=c$ and $b=b_\infty\cos\theta$. And if $h$ is weakly usc, then $c_\infty\geq c$. (This necessary condition is also sufficient.)

Now choose $h_\infty=\pmatrix a_\infty &b_\infty\\
\overline{b_\infty} &c_\infty\endpmatrix$, a positive invertible matrix with $b_\infty\neq 0$, and let $h_n=\pmatrix a_\infty\cos^2\theta &b_\infty\cos\theta\\
\overline{b_\infty}\cos\theta &c_\infty\endpmatrix$ for $n<\infty$. Then $h\in p\asa p$, and $h\geq\ep p$ for some $\ep>0$. A routine calculation shows that $h^{-1}$ is not middle usc. Another such calculation shows that $(h-t_0p)^{-1}$ is not weakly usc for $0<t_0<\ep$. Note that by taking $\theta$ close to 0 and choosing $h_\infty$ appropriately we can arrange that $\sigma(h)$ is contained in a very small interval, say $[1-\del, 1+\del]$. This means that the operator convex functions $f_1(t)=1/t$, $t\in(0,\infty)$, and $f_2(t)=1/(t-t_0)$, $t\in(t_0,\infty)$, do not have better usc properties than claimed in Theorem 4.3, even if they are restricted to small subsets of their natural domains.

We next prove a result about the effect of an operator convex function $f$ on $pM(A)_{\sa}p$, but the hypothesis on the domain of $f$ is much stronger than in Theorem 4.3. The last example shows that the domain hypothesis of 4.3 would not suffice. (Note that since $p\in pM(A)_{\sa}p$, all results are translation invariant.) 
\endexample

\proclaim{Proposition 4.6} Assume $p$ is a closed projection in $A^{**}$, $f$ is an operator convex function on an interval $I$, $\th\in M(A)_\sa$, and $\sigma(\th)\subset I$. Let $h=p\th p$. Then $f(h)$ is middle usc on $p$. If $f(\pi(\th))\leq 0$, where $\pi:M(A)\to M(A)/A$ is the quotient map, then $f(h)$ is strongly usc on $p$.
\endproclaim

\demo{Proof} It is enough to prove the last sentence. We apply 4.2 or 4.3 for the unital $C^*$--algebra $M(A)$, identifying $M(A)^{**}$ with $A^{**}\oplus(M(A)/A)^{**}$. We use the closed projection $\tp=p\oplus\b1_{M(A)/A}$. Then $f(h)\oplus f(\pi(\th))$ is usc on $\tp$, and we claim this implies $f(h)$ is strongly usc on $p$. Thus consider a net $(\var_\al)$ in $\fp$ such that $\var_\al\to \var$ weak$^*$. It is permissible to pass to a subnet such that $\var_\al\to\theta\oplus\psi$ in $F(\tp)$. Since $\theta$ and $\var$ agree on $A$, $\theta=\var$. Now $\var(f(h))+\psi(f(\pi(\th)))\geq\lim\sup \var_\al(f(h))$. Since $f(\pi(\th))\leq 0$, this implies $\var(f(h))\geq\lim\sup\var_\al(f(h))$.
\enddemo

\noindent
{\bf{Note.}} It would be sufficient to assume $I\supset\sigma_{pA^{**}p}(h)\cup\sigma(\pi(\th))$, but this does not lead to a stronger result. If this is the case, 1.2 can be used to find $\th'$ in $M(A)$ such that $\tp\th'\tp=\tp\th\tp$ and $\sigma(\th')\subset I$.

\proclaim{Lemma 4.7} Assume $p$ is a closed projection in $A^{**}$, $A$ is unital, $f$ is a strongly operator convex function on an interval $I$, $h\in p\asa p$, and $\sigma_{pA^{**}p}(h)\subset I$. Then $f(h)$ is usc in $A^{**}$.  Here $f(h)$ is computed in $pA^{**}p$ and is regarded as an element of $A^{**}$ via the inclusion $pA^{**}p \subset A^{**}$.
\endproclaim 

\demo{Proof} We use the integral representation of $f$. Then $f(h)$ can be obtained by substituting $h$ for $x$ in (10), thus obtaining a Bochner integral. Since $((\asa)_{\m})^-$ is a closed cone, it is enough to verify that each value of the integrand is usc in $A^{**}$ and note that $cp$ is usc in $A^{**}$. The former follows from Lemma 4.1(ii).
\enddemo

\proclaim{Theorem 4.8} Assume $p$ is a closed projection in $A^{**}$, $f$ is a strongly operator convex function on an interval $I$, $h\in p\asa p$, and $\sigma_{pA^{**}p}(h)\subset I$.  Define $f(h)$ in $A^{**}$ as in Lemma 4.7.

(i) If $p$ is compact, then $f(h)$ is strongly usc in $A^{**}$.

(ii) If $0\in I$, then $f(h)-f(0)\b1$ is strongly usc in $A^{**}$, whence $f(h)$ is middle usc in $A^{**}$. 

(iii) If $0$ is an endpoint of $I$, then $f(h)$ is weakly usc in $A^{**}$.

\noindent
Conversely, if $f$ is a function which satisfies the conclusion of (i) for all closed faces in the unital case, then $f$ is strongly operator convex.
\endproclaim

\demo{Proof} (i) We apply Lemma 4.7 to $\tA$, identifying $(\tA)^{**}$ with $A^{**}\oplus\bC$, and using the closed projection $p\oplus 0$. Thus $f(h)\oplus 0$ is usc in $(\tA)^{**}$, and this implies $f(h)$ is strongly usc in $A^{**}$.

(ii)  Again we apply Lemma 4.7 to $\tA$, now with the closed projection $p\oplus 1$.  Then $f(h)\oplus f(0)$ is usc in $(\tA)^{**}$.  (Note that $f(0) > 0$ except in the trivial case $f=0$.)  Since $(f(h)-f(0)\b1)\oplus 0 = f(h)\oplus f(0) - f(0)\b1_{(\tA)^{**}}$, this implies $f(h)-f(0)\b1$ is strongly usc in $A^{**}$.

(iii)  This follows from (ii) just as in the proof of Theorem 4.3.

For the converse we use the same algebra $A$, projection $p$, and element $h$ as in the proof of Theorem 4.3.  Then $f(h)$ is the element $x$ of $A^{**}$ given by $x_n = \pmatrix f(a_{\infty}) &0\\
0 &0\endpmatrix$ for $n$ finite and $x_{\infty} = f(h_{\infty})$.   It is easily seen that the upper semicontinuity of $x$ implies (for this $x$) that $x_n\leq x_{\infty}$.  This condition, demanded for all $k,l$, implies the criterion for strong operator convexity given at the beginning of this section.
\enddemo

\proclaim{Theorem 4.9} Assume $p$ is a closed projection in $A^{**}$, $f$ is an operator convex and operator decreasing on an interval $I$, $h$ in $p\asa^{**}p$ is strongly lsc on $p$, and $\sigma_{pA^{**}p}(h)\subset I$. 

(a) (i) If either $p$ is compact or $0\in I$ and $f(0)\leq 0$, then $f(h)$ is strongly usc on $p$.

\quad(ii) If $0\in I$, then $f(h)$  is middle usc on $p$.

\quad(iii) If $0$ is an endpoint of $I$, then $f(h)$ is weakly usc on $p$.

(b) Assume further that $f$ is strongly operator convex.

\quad (i')  If $p$ is compact, then $f(h)$ is strongly usc in $A^{**}$.

\quad (ii') If $0\in I$, then $f(h)$ is middle usc in $A^{**}$.

\quad (iii') If $0$ is an endpoint of $I$, then $f(h)$ is weakly usc in $A^{**}$.

\noindent
There are three other symmetric versions of this result, so that all the cases where $f$ is (strongly) operator convex or (strongly) operator concave and also operator monotone or operator decreasing are covered.
\endproclaim

\demo{Proof} (a) (i) The natural domain of $f$ is an interval of the form $(t_0,\infty)$ or $[t_0,\infty)$, $t_0\leq 0$, so we assume $I$ has this form. Let $h$ be the limit of an increasing net $(h_\al)=(\lam_\al p+pa_\al p)$, as in Proposition 1.3 (ii) and take $\del>0$. By Dini's theorem, $f(h_\al+\del p)$ is defined for $\al$ sufficiently large, and it is the result of applying $f(\del+\lam_\al+\cdot)$ to $pa_\al p$. Assuming, as we may, that $\del+\lam_\al>0$, $f(\del+\lam_\al+0)\leq f(0)\leq 0$ in the non--compact case. In either case, then, Theorem 4.3 implies that $f(h_\al+\del p)$ is strongly usc. Since $f(h_\al+\del p)$ decreases to $f(h+\del p)$, $f(h+\del p)$ is also strongly usc. And since $f(h+\del p)$ converges in norm to $f(h)$, finally $f(h)$ is strongly usc.

(ii) Follows from (i).

(iii) Since $f(h+\del p)$ is obtained by applying $f(\cdot +\del)$ to $h$, part (ii) implies $f(h+\del p)\in(p\Asa p)_\m$, $\forall\del>0$, whence $f(h)\in ((p\Asa p)_{\m})^-$.

(b)(i')  The proof is the same as for the compact case of (a)(i) except that we use Theorem 4.8 instead of Theorem 4.3.

(ii') The proof is the same as for the non-compact case of (a)(ii) except that we use Theorem 4.8 to deduce that $f(\delta p+\lambda_{\alpha} p+pa_{\alpha}p) - f(\delta +\lambda_{\alpha})\b1$ is strongly usc in $A^{**}$.  Since $f(\delta + \lambda_{\alpha}) \le f(0)$, also $f(h_{\alpha} +\delta p) -f(0)\b1$ is strongly usc in $A^{**}$.  Then, as above, it follows that $f(h)-f(0)\b1$ is strongly usc in $A^{**}$.

(iii') follows from (ii') as above.
\enddemo

For the symmetrical versions use one of the functions $\pm f(\pm\cdot)$.

One of the very useful results of \cite{AP} is that for $h\geq\ep\b1, \ep>0$, $h$ is strongly lsc if and only if $h^{-1}$ is weakly usc ([AP,Proposition 3.5]). (Combes \cite{C2} had previously done the unital case, where the complications referred to in the title of \cite{AP} don't occur.) In [B1, Proposition 2.1] we expanded on this by considering also the other two kinds of upper semicontinuity for $h^{-1}$. One part of [B1, Proposition 2.1] is that $h^{-1}$ is never strongly usc unless $A$ is unital. In the case of a closed face the situation is more complicated.
 
\proclaim{Proposition 4.10} Assume $p$ is a closed projection in $A^{**}$, $h\in pA^{**}p$, and $h\geq\ep p$ for some $\ep>0$.  Let $h^{-1}$ denote the inverse of $h$ in $pA^{**}p$, regarded as an element of $A^{**}$.

(i) Then $h^{-1}$ is weakly usc in $A^{**}$ if and only if $h$ is strongly lsc on $p$.

(ii) Then $h^{-1}$ is middle usc in $A^{**}$ if and only if there is $\eta > 0$ such that $h-\eta p$ is strongly lsc on $p$.

(iii) Then $h^{-1}$ is strongly usc in $A^{**}$ if and only if $p$ is compact and $h$ is lsc on $p$.
\endproclaim

\demo{Proof} (i) If $h^{-1}$ is weakly usc in $A^{**}$, then $h^{-1} + \delta \b1$ is weakly usc for all $\delta >0$.  By [AP] $(h^{-1} +\delta \b1)^{-1}$ is strongly lsc in $A^{**}$, whence $p(h^{-1} +\delta \b1)^{-1}p$ is strongly lsc on $p$.  But $p(h^{-1} +\delta \b1)^{-1}p = (h^{-1} +\delta p)^{-1}$ and this converges in norm to $h$ as $\delta \to 0$.  Hence $h$ is strongly lsc on $p$.  The converse follows from Theorem 4.9(iii').

(ii) If $A$ is unital, (ii) follows from (i), so we may assume $A$ non-unital. If $h^{-1}$ is middle usc in $A^{**}$, let $(k_i)_{i \in D}$ be a decreasing net in $\Asa$ which converges to $h^{-1}$.  Let $D' = D \times (0,\infty)$ and define $(i_1,\delta_1) \le (i_2,\delta_2)$ if and only if $i_1 \le i_2$ and $\delta_1 \ge \delta_2$.  Then if $k_{i,\delta} = k_i + \delta \b1$, $(k_{i,\delta})_{(i,\delta) \in D'}$ is also a decreasing net converging to $h^{-1}$.  Then $({k_{i,\delta}}^{-1})$ is an increasing net in $\Asa$.  Since $k_{i,\delta} \ge h^{-1} +\delta \b1$, ${k_{i,\delta}}^{-1} \le (h^{-1} +\delta \b1)^{-1}$ and $p{k_{i,\delta}}^{-1}p \le p(h^{-1} +\delta \b1)^{-1}p = (h^{-1} + \delta p)^{-1} \le h$.  Since for fixed $\delta$, $(p{k_{i,\delta}}^{-1}p)_{i \in D}$ converges to $(h^{-1} +\delta p)^{-1}$, the limit of $(p{k_{i,\delta}}^{-1}p)_{i \in D'}$ must be at least $(h^{-1} +\delta p)^{-1}$, for all $\delta > 0$, and hence the limit is $h$.  Let $\eta > 0$ be the scalar component of ${k_{i_0,\delta_0}}^{-1}$ for some $(i_0,\delta_0)$ in $D'$.  Then $h-\eta p$ is strongly lsc on $p$.  The converse follows from Theorem 4.9(ii') applied to $h-\eta p$ and $f(x)= (x + \eta)^{-1}$.

(iii) If $h^{-1}$ is strongly usc in $A^{**}$, then $h^{-\alpha}$ is strongly lsc in $A^{**}$ for $0 < \alpha < 1$ by [B1, Proposition 2.30(a)].  Since $h^{-\alpha} \to p$ in norm as $\alpha \to 0$, $p$ is strongly usc in $A^{**}$, whence $p$ is compact by [B1, 2.47].  (This argument proves a general result which should have been part of [B1, 2.44(b)] but was neglected by an oversight.)  And the fact that $h$ is lsc on $p$ follows from part (i).  The converse follows from Theorem 4.9(i').
\enddemo

\proclaim{Corollary 4.11} Assume $p$ is a closed projection in $A^{**}$, $h\in p\asa^{**} p$, and $h\geq\ep p$ for some $\ep>0$.

(i) If $h$ is strongly lsc on $p$, then $h^{-1}$ is weakly usc on $p$.

(ii) If $\exists \eta>0$ such that $h-\eta p$ is strongly lsc on $p$, then $h^{-1}$ is middle usc on $p$.
\endproclaim

The fact that the converses don't hold follows from the last sentence of Theorem 4.3, or it can be deduced from Proposition 4.10 by exhibiting $h$ such that $h^{-1}$ is usc on $p$ but not in $A^{**}$. 

It was shown in [B1, Proposition 2.59(a)] that if $f$ is non--linear and operator convex and if both $h$ and $f(h)$ are in $QM(A)_\sa$, then $h$ must be in $M(A)$. Theorem 4.14 below is an analogue with $pAp$ in place of $QM(A)$. (We have found out recently that both Theorem 4.14(i) and [B1, 2.59(a)] can be generalized by allowing $f$ to be merely continuous and strictly convex.  But we don't know of any similar generalization for Theorem 4.14(ii) or [B1, 2.59(b), (c)].)  If $h$ is in $p\asa^{**}p$, $h$ is called {\it{strongly}} {\it{$q$--continuous}} {\it{on}} $p$ if $\chi_F(h)$ is closed whenever $F$ is a closed subset of $\bR$ and $\chi_F(h)$ is compact if in addition $0\not\in F$. Here $\chi_F(h)$ is the spectral projection, computed in $pA^{**}p$. Let $SQC(p)$ denote the complex span of the elements strongly $q$--continuous on $p$. If was shown in [B1, Theorem 3.43 (a)] that $SQC(p)=\{pa: a\in A\;\text{and}\;pa=ap\}$, whence $SQC(p)$ is a $C^*$--algebra. It was shown in the proof that 4 implies 1 in [B3, Theorem 3.1] that if $h$ and $h^2$ are both in $p\asa p$, then $h\in SQC(p)$. Thus $SQC(p)$ is the largest $C^*$--algebra contained in $pAp$, just as $M(A)$ is the largest $C^*$--algebra contained in $QM(A)$. (The last fact follows from [AP, Proposition 4.4].)

\proclaim{Lemma 4.12} Assume $p$ is a closed projection in $A^{**}$, $A$ is unital, and $h\in p\asa p$.

(i) If $h^2\in p\asa p$, then $h\in SQC(p)$.

(ii) If $h\geq\ep p$ for some $\ep>0$ and $h^{-1}\in pAp$, then $h\in SQC(p)$.
\endproclaim

\demo{Proof} (i) Follows from \cite{B3} as just noted.

(ii) Two methods can be used to obtain $h^{1\over 2}$: Apply the operator concave function, $f_1(t)=t^{1\over 2}$, to $h$ or apply the operator convex function, $f_2(t)=t^{-{1\over 2}}$ to $h^{-1}$. By 4.2  or 4.3 we see that $h^{1\over 2}$ is both lsc and usc on $p$, whence $h^{1\over 2}\in p\asa p$. By (i), $h^{1\over 2}\in SQC(p)$, and thus also $h\in SQC(p)$.
\enddemo

\proclaim{Lemma 4.13} Assume $p$ is a closed projection in $A^{**}$, $f$ is non--linear on an interval $I$, $h\in p\asa^{**}p$, $\sigma_{pA^{**}p}(h)\subset I$, and $A$ is unital.

(i) If $f$ is operator convex and both $h$ and $f(h)$ are in $p\asa p$, then $h\in SQC(p)$.

(ii) If $f$ is operator convex and operator decreasing, $h$ is lsc on $p$, and $f(h)\in pAp$, then $h\in SQC(p)$.
\endproclaim

\demo{Proof} (i) We can obtain $f(h)$ by substituting $h$ for $x$ in formula (9). A basic principle is that if both $k_1$ and $k_2$ are usc on $p$ and $k_1+k_2$ is in $pAp$, then $k_1$ and $k_2$ are both in $pAp$ (since they are then lsc on $p$). We can obtain this situation, using 4.2 or 4.3 by writing $f$ as the sum of two operator convex functions; and the latter can be achieved by partioning $\bR\setminus I$ into two Borel sets. The conclusion is that for $t_0$ in the closed support of $\mu_+$ or $\mu_-$, and for a neighborhood $E$ of $t_0$, 
$${1\over \mu_{\pm}(E)}\int_E{\pm(h-x_0)^2\over (tp-h)(t-x_0)^2}d\mu_{\pm}(t)\in p\asa p.$$
Letting $E$ shrink to $\{t_0\}$, we see that ${(h-x_0)^2\over (t_0 p-h)(t_0-x_0)^2}$ is in $pAp$, whence ${1\over t_0 p-h}$ is in $pAp$. Then 4.12 (ii) implies $h$ is in $SQC(p)$. If both $\mu_+$ and $\mu_-$ are 0, then $a>0$, and we obtain $h^2\in pAp$ and cite 4.12 (i). 

(ii) We proceed in a similar way, using (11) (for $-f$) instead of (9) and 4.9 (i) instead of 4.3. The conclusion is (recall that $I=(s,\infty)$ or $[s,\infty)$) that for some $t_0$, $h-t_0p\geq\ep p, \ep>0$, and $(h-t_0p)^{-1}\in pAp$. But then 4.2 or 4.3 applied to $(h-t_0p)^{-1}$ implies that $h-t_0p$ is usc on $p$, whence $h\in pAp$. Then part (i) applies.
\enddemo

\proclaim{Theorem 4.14} Assume $p$ is a closed projection in $A^{**}$, $f$ is a continuous function on an interval $I$, $h\in p\asa^{**}p$, and $\sigma_{pA^{**}p}(h)\subset I$. Also assume either $p$ is compact or $0\in I$ and $f(0)=0$.

(i) If $f$ is non--linear and operator convex and if both $h$ and $f(h)$ are in $pAp$, then $h\in SQC(p)$.

(ii) If $f$ is non--linear and operator convex and operator decreasing, $h$ is strongly lsc on $p$, and $f(h)$ is in $pAp$, then $h\in SQC(p)$.

\noindent
As in Theorem 4.9, there are three other symmetrical versions of part (ii).
\endproclaim

\demo{Proof} We apply Lemma 4.13 to $\tA$, identifying $\tA^{**}$ with $A^{**}\oplus\bC$. If $p$ is compact, then $p\oplus 0$ is closed in $\tA^{**}$, and we conclude that $h\oplus 0$ is in $SQC(p\oplus 0)$, whence $h$ is in $SQC(p)$. Otherwise $\tp=p\oplus 1$ is closed in $\tA^{**}$, and, using the fact that $f(h\oplus 0)=f(h)\oplus 0$, we conclude that $h\oplus 0$ is in $SQC(p\oplus 1)$. If $F$ is closed in $\bR$ and $0\in F$, this yields the fact that $\chi_F(h)\oplus 1$ is closed in $\tA^{**}$, which implies that $\chi_F(h)$ is closed in $A^{**}$. If $0\not\in F$, then $\chi_F(h)\oplus 0$ is closed in $\tA^{**}$, which implies that $\chi_F(h)$ is compact in $A^{**}$.                              \enddemo

\example{Remark 4.15} Assume $p$ is not compact. Is the hypothesis that $0\in I$ and $f(0)=0$ really necessary? It is not hard to see that it is impossible to have $f(0)<0$, given the other hypotheses. Also if either $0\notin I$ or $f(0)>0$, it can be seen that the other hypotheses of 4.14 (i) imply that $p$ satisfies condition (1) of \S1. Further analysis of Example 1.8, where of course (1) is satisfied, yields counterexamples for both parts where $0\notin I$ and counterexamples where $f(0)>0$. Also note that when $p$ is not compact, it is impossible to have $h \in SQC(p)$ and $0 \notin \sigma (h)$, and it is also impossible to have $h \in SQC(p)$, $f(h) \in pAp$, and $f(0) \ne 0$.  
\endexample

We have given several examples in the attempt to provide evidence for the sharpness of our results; but, although in many cases we have shown that hypotheses can't simply be dropped from the theorems, we do not know whether the results of this section are as sharp as the corresponding results in \cite{B1}, as documented on pages 899--902 of \cite{B1}.

\S 5. Final remarks and questions.

In [B1, Propostion 2.2 and Theorem 3.27] we proved that if weak and middle semicontinuity coincide for elements of $A^{**}_{\text {sa}}$, then several other properties hold.  One of these is that every positive element $(\Asa)^{\text {m}}$ is strongly lsc.  The same is not true for semicontinuity on $p$, as shown by the closed face of Example 1.8.  On the other hand, Proposition 4.10 can be used together with the arguments in [B1, Proposition 2.2] to show that if $(\Asa)^{\text {m}} = ((\Asa)^{\text {m}})^-$, then $(p\Asa p)^{\text {m}} = ((p\Asa p)^{\text {m}})^-$ for all closed projections $p$ in $A^{**}$, and every positive element of $(p\Asa p)^{\text {m}}$ is strongly lsc on $p$.  

Question 1.  Are there any other special properties that hold when, for a particular closed face, $(p\Asa p)^{\text {m}} = ((p\Asa p)^{\text {m}})^-$?

In [B1] we sometimes considered the following partition of the class of $C^*$-algebras:

1.  Unital $C^*$-algebras.

2.  Non-unital $C^*$-algebras for which $\Asa^{\text {m}} = (\Asa^{\text {m}})^-$.

3.  Non-unital $C^*$-algebras for which $\Asa^{\text {m}} \ne (\Asa^{\text {m}})^-$.

This trichotomy is appropriate for semicontinuity theory.  In this paper, despite the fact that there is a tremendous variety within the class of closed faces of $C^*$-algebras, our theorems distinguished only between compact projections and non-compact closed projections.

Question 2.  What are the "right" subclasses of the class of closed faces of $C^*$-algebras to consider in connection with semicontinuity theory?

\bye